\documentclass[11pt]{article}
\usepackage[mathscr]{eucal}
\DeclareMathAlphabet{\mathpzc}{OT1}{pzc}{m}{it}
\usepackage{bbm}
\usepackage{dsfont}
\usepackage{color}
\usepackage{mathrsfs}
\usepackage{amssymb}
\usepackage{calligra}
\usepackage{fontenc}
\usepackage{amsbsy}
\newcommand{\bfscr}[1]{{\pmb{\mathscr{#1}}}}
\newcommand{\bfpzc}[1]{{\pmb{\mathpzc{#1}}}}
\newcommand{\bfsf}[1]{{\textbf{\textsf{#1}}}}
\newcommand{\sub}[1]{{\mbox{\footnotesize $#1$}}}
\usepackage{graphicx}
\graphicspath{%
    {converted_graphics/}
    {/}
    {C:/Users/GALDI/Desktop/Vibration/}
}
\begin{document}

\newcommand{\sfW}{{\sf W}}
\newcommand{\sfD}{{\sf D}}
\newcommand{\sfL}{{\sf L}}
\newcommand{\sfw}{{\sf w}}
\newcommand{\sfp}{{\sf p}}
\newcommand{\simge}{\ba{cc}\vspace*{-2.4mm}>\\ \sim\ea }
\newcommand{\simle}{\ba{cc}\vspace*{-2.4mm}<\\ \sim\ea }
\newcommand{\Cdot}{\!\cdot\!}
\newcommand{\sq}{{$\sqcap\!\!\!\!\sqcup$}}
\newcommand{\Eu}{{\rm I\,\!\! E}}
\newcommand{\Io}{\Int{\Omega}{}}
\newcommand{\Id}{\Int{\cald}{}}
\newcommand{\Div}{\mbox{\rm div}\,}
\newcommand{\tr}{\mbox{\rm tr}\,}
\newcommand{\grad}{\mbox{\rm grad}\,}
\newcommand{\supp}{\mbox{\rm supp}\,}
\newcommand{\curl}{\mbox{\rm curl}\,}
\newcommand{\Ido}{\Int{\partial\Omega}{}}
\newcommand{\IdS}{\Int{\Sigma}{}}
\newcommand{\Oint}[2]{{\displaystyle \oint_{#1}^{#2}}}
\newcommand{\Int}[2]{{\displaystyle \int_{ #1}^{ #2}}}
\newcommand{\Lim}[1]{{\displaystyle \lim_{ #1}}}
\newcommand{\Limsup}[1]{{\displaystyle \limsup_{\footnotesize #1}}}
\newcommand{\Liminf}[1]{{\displaystyle \liminf_{\footnotesize #1}}}
\newcommand{\Sup}[1]{{\displaystyle \sup_{#1}}}
\newcommand{\Inf}[1]{{\displaystyle \inf_{#1}}}
\newcommand{\Max}[1]{{\displaystyle \max_{#1}}}
\newcommand{\Min}[1]{{\displaystyle \min_{#1}}}
\newcommand{\Sum}[2]{{\displaystyle \sum_{#1}^{#2}}}
\newcommand{\Prod}[2]{{\displaystyle \prod_{#1}^{#2}}}
\newcommand{\BCup}[2]{{\displaystyle \bigcup_{#1}^{#2}}}
\newcommand{\BCap}[2]{{\displaystyle \bigcap_{#1}^{#2}}}
\newcommand{\Frac}[2]{\displaystyle{\frac{\displaystyle{#1}}{\displaystyle{#2}}}}
\newcommand{\norm}[1]{\left\|{#1}\right\|}
\newcommand{\Norm}[1]{\langle\langle{#1}\rangle\rangle_q}
\newcommand{\No}[1]{\langle\!\langle{#1}\rangle\!\rangle}
\newcommand{\NO}[1]{{\langle{#1}\rangle}_{\lambda,q}}
\newcommand{\beea}{\begin{eqnarray}}
\newcommand{\eeea}{\end{eqnarray}}
\newcommand{\ms}{\medskip\smallskip}
\newcommand{\bs}{\bigskip}
\newcommand{\ps}{\par\smallskip}
\newcommand{\bfe}{{\mbox{\boldmath $e$}} }
\newcommand{\pni}{{\par\noindent}}
\newcommand{\bfq}{{\mbox{\boldmath $q$}} }
\newcommand{\bfz}{{\mbox{\boldmath $z$}} }
\newcommand{\0}{{\mbox{\boldmath $0$}} }
\newcommand{\LE}{\!\!\!&\le&\!\!\!}
\newcommand{\BL}[1]{{\par\smallskip{\bf Lemma #1.}}}
\newcommand{\BT}[1]{{\par\smallskip{\bf Theorem #1.}}}
\newcommand{\Ln}{[\!|}
\newcommand{\Rn}{|\!]}
\newcommand{\C}{{\sf c}}
\newcommand{\n}[1]{{\Ln{#1}\Rn}} 
\newcommand{\nq}[1]{{\Ln{#1}\Rn}_{q}} 
\newcommand{\nqr}[1]{{\Ln{#1}\Rn}_{q,r}} 
\newcommand{\Nq}[1]{{\langle{#1}\rangle}_{q}} 
\newcommand{\Nql}[1]{{\langle{#1}\rangle}_{\lambda,q}} 
\newcommand{\Nqr}[1]{{\langle{#1}\rangle}_{q,r}}
\newcommand{\N}[1]{{|\!\!|\!\!|\,{#1}\,|\!|\!\!|_2}}
\newcommand{\EA}[2]{$$#1$$%
\vspace{-6.mm}
\begin{equation}
\end{equation}
\vspace{-6.mm}
$$
#2
\setlength{\belowdisplayskip}{3mm}
\setlength{\belowdisplayshortskip}{3mm}
$$
}
\newcommand{\A}[2]{$$#1$$%
\vspace{-4.mm}
$$
#2
\setlength{\belowdisplayskip}{3mm}
\setlength{\belowdisplayshortskip}{3mm}
$$
}
\newcommand{\BF}{\begin{footnotesize}}
\newcommand{\EF}{\end{footnotesize}}
\setlength{\jot}{.15in}
\newcommand{\pde}[2]{{\displaystyle \frac{\mbox{$\partial #1$}}{\mbox{$\partial #2$}}}}
\newcommand{\ode}[2]{{\displaystyle \frac{\mbox{$d #1$}}{\mbox{$d #2$}}}}
\newcommand{\f}[2]{\frac{\mbox{$#1$}}{\mbox{$ #2$}}}
\newcommand{\bi}{\begin{itemize}}
\newcommand{\ei}{\end{itemize}}
\newcommand{\ed}{\end{document}}
\newcommand{\be}{\begin{equation}}
\newcommand{\ba}{\begin{array}}
\newcommand{\ea}{\end{array}}
\newcommand{\ee}{\end{equation}}
\newcommand{\eeq}[1]{\label{eq:#1}\end{equation}}
\newcommand{\real}{{\mathbb R}}
\newcommand{\compl}{{\mathbb C}}
\def\Id{\mbox{\boldmath $1$}}
\def\zero{\mbox{\boldmath $0$}}
\newcommand{\PP}{{\rm I\!\!\,P}}
\newcommand{\nat}{{\mathbb N}}
\newcommand{\bfpsi}{\mbox{\boldmath $\psi$}}
\newcommand{\bfchi}{\mbox{\boldmath $\chi$}}
\newcommand{\bfomega}{\mbox{\boldmath $\omega$}}
\newcommand{\bfvaromega}{\mbox{\boldmath $\varpi$}}
\newcommand{\bfOmega}{\mbox{\boldmath $\Omega$}}
\newcommand{\bfTheta}{\mbox{\boldmath $\Theta$}}
\newcommand{\bfxi}{\mbox{\boldmath $\xi$}}
\newcommand{\bfmu}{\mbox{\boldmath $\mu$}}
\newcommand{\bfx}{\mbox{\boldmath $x$}}
\newcommand{\bfy}{\mbox{\boldmath $y$}}
\newcommand{\bfPsi}{\mbox{\boldmath $\Psi$}}
\newcommand{\bfphi}{\mbox{\boldmath $\varphi$}}
\newcommand{\bfhi}{\mbox{\boldmath $\phi$}}
\newcommand{\bfPhi}{\mbox{\boldmath $\Phi$}}
\newcommand{\bfv}{{\mbox{\boldmath $v$}} }
\newcommand{\bfu}{{\mbox{\boldmath $u$}} }
\newcommand{\bfuf}{{\mbox{\footnotesize\boldmath $u$}} }
\newcommand{\bfw}{{\mbox{\boldmath $w$}} }
\newcommand{\bff}{{\mbox{\boldmath $f$}} }
\newcommand{\bfa}{{\mbox{\boldmath $a$}} }
\newcommand{\bfi}{{\mbox{\boldmath $i$}} }
\newcommand{\bfj}{{\mbox{\boldmath $j$}} }
\newcommand{\bfc}{{\mbox{\boldmath $c$}} }
\newcommand{\bfo}{{\mbox{\boldmath $o$}} }
\newcommand{\bfp}{{\mbox{\boldmath $p$}} }
\newcommand{\bfkp}{{\mbox{\footnotesize{\boldmath $k$}}} }
\newcommand{\bfka}{{\mbox{\footnotesize{\boldmath $k^*$}}} }
\newcommand{\bft}{{\mbox{\boldmath $t$}} }
\newcommand{\bfd}{{\mbox{\boldmath $d$}} }
\newcommand{\bfl}{{\mbox{\boldmath $l$}} }
\newcommand{\bfr}{{\mbox{\boldmath $r$}} }
\newcommand{\bfk}{{\mbox{\boldmath $k$}} }
\newcommand{\bfA}{{\mbox{\boldmath $A$}} }
\newcommand{\bfS}{{\mbox{\boldmath $S$}} }
\newcommand{\bfO}{{\mbox{\boldmath $O$}} }
\newcommand{\bfM}{{\mbox{\boldmath $M$}} }
\newcommand{\bfP}{{\mbox{\boldmath $P$}} }
\newcommand{\bfB}{{\mbox{\boldmath $B$}} }
\newcommand{\bfR}{{\mbox{\boldmath $R$}} }
\newcommand{\bfC}{{\mbox{\boldmath $C$}} }
\newcommand{\bfD}{{\mbox{\boldmath $D$}} }
\newcommand{\bfQ}{{\mbox{\boldmath $Q$}} }
\newcommand{\bfZ}{{\mbox{\boldmath $Z$}} }
\newcommand{\bfG}{{\mbox{\boldmath $G$}} }
\newcommand{\bfE}{{\mbox{\boldmath $E$}} }
\newcommand{\bfX}{{\mbox{\boldmath $X$}} }
\newcommand{\bfY}{{\mbox{\boldmath $Y$}} }
\newcommand{\bfH}{{\mbox{\boldmath $H$}} }
\newcommand{\bfI}{{\mbox{\boldmath $I$}} }
\newcommand{\bfJ}{{\mbox{\boldmath $J$}} }
\newcommand{\bfN}{{\mbox{\boldmath $N$}} }
\newcommand{\bfh}{{\mbox{\boldmath $h$}} }
\newcommand{\bfm}{{\mbox{\boldmath $m$}} }
\newcommand{\bfone}{{\mbox{\boldmath $1$}} }
\newcommand{\hs}{{\rm I}\!\!\,{\rm R}^3_+}
\newcommand{\cala}{{\cal A}}
\newcommand{\calb}{{\cal B}}
\newcommand{\calc}{{\cal C}}
\newcommand{\cald}{{\cal D}}
\newcommand{\cale}{{\cal E}}
\newcommand{\calf}{{\cal F}}
\newcommand{\calg}{{\cal G}}
\newcommand{\calh}{{\cal H}}
\newcommand{\cali}{{\cal I}}
\newcommand{\calj}{{\cal J}}
\newcommand{\calk}{{\cal K}}
\newcommand{\call}{{\cal L}}
\newcommand{\calm}{{\cal M}}
\newcommand{\caln}{{\cal N}}
\newcommand{\calo}{{\cal O}}
\newcommand{\calp}{{\cal P}}
\newcommand{\calq}{{\cal Q}}
\newcommand{\calr}{{\cal R}}
\newcommand{\cals}{{\cal S}}
\newcommand{\calt}{{\cal T}}
\newcommand{\calu}{{\cal U}}
\newcommand{\calv}{{\cal V}}
\newcommand{\calx}{{\cal X}}
\newcommand{\caly}{{\cal Y}}
\newcommand{\calw}{{\cal W}}
\newcommand{\calz}{{\cal Z}}
\newcommand{\bfsigma}{\mbox{\boldmath $\sigma$}}
\newcommand{\bfSigma}{\mbox{\boldmath $\Sigma$}}
\newcommand{\bftau}{\mbox{\boldmath $\tau$}}
\newcommand{\bfeta}{\mbox{\boldmath $\eta$}}
\newcommand{\bfT}{{\mbox{\boldmath $T$}} }
\newcommand{\bfV}{{\mbox{\boldmath $V$}} }
\newcommand{\bfU}{{\mbox{\boldmath $U$}} }
\newcommand{\bfW}{{\mbox{\boldmath $W$}} }
\newcommand{\bfF}{{\mbox{\boldmath $F$}} }
\newcommand{\bfK}{{\mbox{\boldmath $K$}} }
\newcommand{\bfL}{{\mbox{\boldmath $L$}} }
\newcommand{\bfb}{{\mbox{\boldmath $b$}} }
\newcommand{\bfg}{{\mbox{\boldmath $g$}} }
\newcommand{\bfn}{{\mbox{\boldmath $n$}} }
\newcommand{\bfs}{{\mbox{\boldmath $s$}} }
\newcommand{\cf}{{\it cf.} }
\newcommand{\io}{\int_\Omega}
\newcommand{\1}{\item[({\it i})]}
\newcommand{\2}{\item[({\it ii})]}
\newcommand{\3}{\item[({\it iii})]}
\newcommand{\4}{\item[({\it iv})]}
\newcommand{\5}{\item[({\it v})]}
\newcommand{\6}{\item[({\it vi})]}
\newcommand{\7}{\item[({\it vii})]}
\newcommand{\8}{\item[({\it viii})]}
\newcommand{\9}{\item[({\it xi})]}
\newcommand{\ido}{\int_{\partial\Omega}}
\newcommand{\half}{\mbox{$\frac{1}{2}$}}
\def\parallel{\|}
\def\mid{|}
\def\Bbb R{\real}
\def\hat{\widehat}
\def\tilde{\widetilde}
\def\bar{\overline}
\newcommand{\threehalves}{3\over 2}
\newcommand{\bfPi}{\mbox{\boldmath $\Pi$}}
\newcommand{\bfXi}{\mbox{\boldmath $\Xi$}}
\newcommand{\bfalpha}{\mbox{\boldmath $\alpha$}}
\newcommand{\bfbeta}{\mbox{\boldmath $\beta$}}
\newcommand{\bfgamma}{\mbox{\boldmath $\gamma$}}
\newcommand{\bfdelta}{\mbox{\boldmath $\delta$}}
\newcommand{\bfzeta}{\mbox{\boldmath $\zeta$}}
\newcommand{\bfUpsilon}{\mbox{\boldmath $\Upsilon$}}
\newcommand{\bfGamma}{\mbox{\boldmath $\Gamma$}}
\newcommand{\bfcala}{\mbox{\boldmath ${\cal A}$}}
\newcommand{\bfcalm}{\mbox{\boldmath ${\cal M}$}}
\newcommand{\bfcaln}{\mbox{\boldmath ${\cal N}$}}
\newcommand{\bfcalq}{\mbox{\boldmath ${\cal Q}$}}
\newcommand{\bfcalb}{\mbox{\boldmath ${\cal B}$}}
\newcommand{\bfcalc}{\mbox{\boldmath ${\cal C}$}}
\newcommand{\bfcali}{\mbox{\boldmath ${\cal I}$}}
\newcommand{\bfcalg}{\mbox{\boldmath ${\cal G}$}}
\newcommand{\bfcalh}{\mbox{\boldmath ${\cal H}$}}
\newcommand{\bfcalk}{\mbox{\boldmath ${\cal K}$}}
\newcommand{\bfcalt}{\mbox{\boldmath ${\cal T}$}}
\newcommand{\bfcalx}{\mbox{\boldmath ${\cal X}$}}
\newcommand{\bfcall}{\mbox{\boldmath ${\cal L}$}}
\newcommand{\bfcalf}{\mbox{\boldmath ${\cal F}$}}
\newcommand{\bfcalr}{\mbox{\boldmath ${\cal R}$}}
\newcommand{\bfcals}{\mbox{\boldmath ${\cal S}$}}
\newcommand{\bfcalw}{\mbox{\boldmath ${\cal W}$}}
\newcommand{\bfcalu}{\mbox{\boldmath ${\cal U}$}}
\newcommand{\bfcalv}{\mbox{\boldmath ${\cal V}$}}
\newcommand{\bfcalz}{\mbox{\boldmath ${\cal Z}$}}
\pagenumbering{roman}
\newcommand{\art}[6]{{\I[{\sc #1,}] {#2}, {\it #3}, {\bf #4}, {#5} {[#6]}}}
\newcommand{\ED}{\end{description}}
\newcommand{\I}{\item }
\newcommand{\ra}{\rm a}
\newcommand{\rb}{\rm b}
\newcommand{\rc}{\rm c}
\newcommand{\s}{{\sf s}}
\newcommand{\Hsp}{{\rm I}\!\!\,{\rm R}^n_+}
\newcommand{\Hsn}{{\rm I}\!\!\,{\rm R}^n_-}
\newcommand{\po}[1]{\mbox{$\displaystyle \frac{\mbox{$\partial #1$}}
{\mbox{$\partial x_{1}$}}$}}
\newcommand{\PO}[1]{\mbox{$\displaystyle \frac{\mbox{$\partial #1$}}
{\mbox{$\partial y_{1}$}}$}}
\newcommand{\OP}{\left(\Delta+2\lambda\PO{}\right)}
\newcommand{\op}{\left(\Delta+2\lambda\po{}\right)}
\newcommand{\ft}[1]{
\Frac{1}{(2\pi)^{n/2}}\Int{{\Bbb R}^{n}}{}e^{i{\bf x}\cdot \bfxi}
#1(\xi)d\xi}
\newcommand{\Ft}[1]{
\Frac{1}{2\pi}\Int{{\Bbb R}^{2}}{}e^{i{x}\cdot \xi}
#1(\xi)d\xi}
\newcommand{\Z}{\item[({\it a})]}
\newcommand{\B}{\item[({\it b})]}
\newcommand{\D}{\item[({\it d})]}
\newcommand{\E}{\item[({\it e})]}
\newcommand{\G}{\item[({\it g})]}
\newcommand{\Š}{\`e}
\newcommand{\…}{\`a}
\newcommand{\•}{\`o}
\newcommand{\—}{\`u}
\newcommand{\}{\`{\i}}
\def\tag{\renewcommand{\theequation}}
\newcommand{\Footnote}{~\footnote}
\newcommand{\ie}{{\it i.e.}}
\newcommand{\dist}{\mbox{\rm dist\,}}
\newcommand{\const}{\mbox{\rm const}}
\newcommand{\trace}{\mbox{\rm trace}}
\newcommand{\Bo}{\par\hfill{$\Box$}\par\noindent}
\newcommand{\Nor}[1]{\langle{#1}\rangle_q}
\newcommand{\vs}{\vspace*{.5cm}\par\noindent}
\newcommand{\Vs}{\vspace*{.6cm}\par\noindent}
\newcommand{\Vvs}{\vspace*{.7cm}\par\noindent}
\newcommand{\VVs}{\vspace*{.8cm}\par\noindent}
\newtheorem{definition}{Definition}[section]
\newcommand{\Bd}{\begin{definition}\begin{rm}}
\newcommand{\Ed}{\end{rm}\end{definition}}
\newtheorem{remark}{Remark}[section]
\newcommand{\Br}{\begin{remark}\begin{rm}}
\newcommand{\Er}{\end{rm}\end{remark}}
\newtheorem{proposition}{Proposition}[section]
\newcommand{\Bp}{\begin{proposition}\begin{sl}}
\newcommand{\EP}[1]{\end{sl}\label{proposition:#1}\end{proposition}}
\newcommand{\propref}[1]{{\rm Proposition \ref{proposition:#1}}}
\newcommand{\Bt}{\begin{theorem}\begin{sl}}
\newcommand{\Et}{\end{sl}\end{theorem}}
\newcommand{\Bl}{\begin{lemma}\begin{sl}}
\newcommand{\El}{\end{sl}\end{lemma}}
\newtheorem{theorem}{Theorem}[section]
\newtheorem{lemma}{Lemma}[section]
\newtheorem{corollary}{Corollary}[section]
\newcommand{\eqref}[1]{{\rm (\ref{eq:#1})}}
\newcommand{\Bc}{\begin{corollary}\begin{sl}}
\newcommand{\Ec}{\end{sl}\end{corollary}}
\newcommand{\ET}[1]{\end{sl}\label{theorem:#1}\end{theorem}}
\newcommand{\EDD}[1]{\end{rm}\label{definition:#1}\end{definition}}
\newcommand{\EL}[1]{\end{sl}\label{lemma:#1}\end{lemma}}
\newcommand{\theoref}[1]{{\rm Theorem \ref{theorem:#1}}}
\newcommand{\defref}[1]{{\rm Definition \ref{definition:#1}}}
\newcommand{\ER}[1]{\end{rm}\label{remark:#1}\end{remark}}
\newcommand{\EC}[1]{\end{sl}\label{corollary:#1}\end{corollary}}
\newcommand{\remref}[1]{{\rm Remark \ref{remark:#1}}}
\newcommand{\cororef}[1]{{\rm Corollary \ref{corollary:#1}}}
\newcommand{\lemmref}[1]{{\rm Lemma \ref{lemma:#1}}}
\newcommand{\essup}[1]{{\rm ess}\,{{\displaystyle \sup_{\hspace*{-5mm}{#1}}}}}

\renewcommand{\i}{{\rm i}}

\pagenumbering{arabic}
\newcommand{\QED}{{\par\par\hfill$\square$\par}}
\renewcommand{\thefootnote}{(\arabic{footnote})}
\title{On Time-Periodic Bifurcation of a Sphere Moving under Gravity\\ in a Navier-Stokes Liquid}
\author{Giovanni P Galdi \smallskip \\ { \small Department of Mechanical Engineering and Materials Science}\\ { \small University of Pittsburgh, USA}}
\date{}
\maketitle
\begin{abstract}We provide sufficient conditions for the occurrence of time-periodic Hopf bifurcation for the coupled system constituted by a rigid sphere, $\mathscr S$, freely moving under gravity in a Navier-Stokes liquid. Since the region of flow is unbounded (namely, the whole space outside $\mathscr S$), the main difficulty consists in finding the appropriate functional setting where general theory may apply. In this regard, we are able to show that the problem can be formulated as a suitable system of coupled operator equations in Banach spaces, where the relevant operators are Fredholm of index 0. In such a way, we can use the theory recently introduced by the author, and give sufficient conditions for time-periodic bifurcation to take place.   
\end{abstract}
\noindent
{\bf Keywords.} Fluid-Structure Interaction; Navier-Stokes Equations; Hopf Bifurcation; Falling Sphere.
\medskip\par\noindent
{\bf AMS Subject Classification.} 35Q30, 76D05, 35B32, 76T99.
\renewcommand{\theequation}{\arabic{section}.\arabic{equation}}
\setcounter{section}{0}
\section*{Introduction} The motion of spheres falling or rising in a viscous liquid has long been recognized as a fundamental topic of research, not only for its intrinsic interest but also for its role in applied sciences; see  \cite{Jenny1,Nak,Rao,Scog,ST,Tachi,Tan} and the bibliography therein. Even though the dynamics may be different depending on whether the sphere is {\em light} (ascending) or {\em heavy} (falling) \cite{Kara}, its {\em qualitative} behavior is rather similar in both cases. More specifically, let $\rho_{\mathscr S}$  and $\rho_{\mathscr L}$ be the density of the sphere and of the liquid, respectively, and denote by $\lambda$ a suitable non-dimensional number depending on $|\rho_{\mathscr S}/\rho_{\mathscr L}-1|$ ({\em Galilei number}); see \eqref{Gal}. Then, experimental and numerical tests  \cite{Jenny1,Nak,Tachi,Tan} show that steady regimes occur as long as $\lambda$ is not too ``large." Precisely, in a first range of Galilei numbers, the sphere merely translates (no spin) with  constant translational velocity, $\bftau_0$, parallel to the direction of gravity, $\bfe$. In this situation, the liquid flow is axisymmetric around $\bfe$. For  $\lambda$ above a first critical value, there is a breaking of symmetry from axisymmetry to planar symmetry. However, the motion is still translatory but now $\bftau_0$ and $\bfe$ are no longer parallel. If $\lambda$ is further increased to some higher critical value, $\lambda_{\sf c}$ (say), then a second (Hopf)       bifurcation occurs: the flow  ceases to be steady and becomes instead purely  time-periodic,  with the constant translatory motion of the sphere giving way to an oscillating
oblique movement. At even greater values of $\lambda$, a chaotic regime eventually sets in.  
\par
The objective of this paper is to furnish a rigorous mathematical contribution to the interpretation of some of the above bifurcation phenomena.
\par
We recall that, from a strict mathematical viewpoint, the study of  bifurcation, in fluid mechanics as in other branches of mathematical physics, presents a fundamental challenge. It consists in determining the {\em appropriate} functional setting where the problem can be formulated in order to be addressed by  general abstract theory. In the situation at hand, this aspect is particularly intriguing, since the region of flow is {\em unbounded in all directions}, which implies that 0 is a point of the {\em essential spectrum} of the relevant linearized operator \cite{FaNe1}. This is a crucial and well--known problem \cite{Bab,Bab1}   that prevents one from  using  classical approaches that are, instead, very successful in the case of {\em bounded} flow, where the above spectral issue is absent  \cite{Satt, HaIo} \cite[Secs. 72.7--72.9]{Z1}.\par  
In the past few years, we have introduced a new, general approach to bifurcation that allows us to overcome the above problem and to provide sufficient (and necessary) conditions for the occurrence of bifurcation, in both steady-state and time-periodic cases, and for both bounded and unbounded flow \cite{Gafur,GaBif}. The basic idea of this approach is to formulate the problem not in classical Sobolev spaces but, instead, in {\em homogeneous} Sobolev spaces, characterized by the property that the various derivatives involved may have {\em different} summability properties  in the neighborhood of spatial infinity. In such a framework, the spectral issue mentioned earlier on is totally absent. In \cite{Gace,GaBif,GaArma} we have  employed this approach to study  bifurcation properties of a Navier--Stokes liquid past a {\em fixed} body, namely, when the body is kept in a given  configuration by suitable forces and torques. More recently, we have used the new method to study steady-state bifurcation of a falling (or ascending) sphere in a viscous liquid  under the action of gravity \cite{GaJMP}, which is a {\em bona fide} fluid-structure interaction problem. In particular, we have shown that a requirement for the occurrence of the above type of bifurcation with $\bftau_0$ parallel to $\bfe$, is that the relevant linearized operator, defined in an {\em appropriate} function space, has 1 as a simple eigenvalue crossing the imaginary axis at ``non-zero speed" (transversality condition).   Remarkably, this requirement {\em formally} coincides with the classical generic bifurcation condition for a flow in a {\em bounded} domain \cite{Satt}.     
\par
In this article we continue and --to an extent-- complete the research initiated in \cite{GaJMP}, by investigating the occurrence of time-periodic bifurcation of the { coupled} system sphere-liquid under the action of gravity. More precisely, we show that, once the problem is formulated in the {\em appropriate functional setting},  we can obtain a bifurcation criterion along the lines  of the classical Hopf theory. Namely, it suffices that the relevant linearized operator has a non-resonating, simple imaginary eigenvalue satisfying the transversality condition. In order to reach this goal, we follow \cite{GaBif,GaArma} and split the  unknowns into the sum of their average over a period plus an oscillatory component. In this way, the original problem transforms into a coupled system of nonlinear elliptic-parabolic equations. We then prove that such a system can be written as two coupled operator equations in suitable spaces, with the relevant operators satisfying all the assumptions of the abstract theory introduced in \cite{GaBif}, which thus provides the desired bifurcation results.
\par
The plan of the paper is as follows. In Section 1 we give the precise formulation of the problem and collect some standard notation. In Section 2, we present  the abstract time-periodic bifurcation result proved in \cite{GaBif}; see \theoref{3.1_ar}. Successively, in Section 3, we recollect some fundamental function spaces introduced in \cite{Gah,ALS} and, for some of them, recall their relevant properties. The following Section 4 is dedicated to the existence of a unique steady-state solution branch parametrized in the Galilei number $\lambda$. To this end, we first show, in \theoref{2.1}, that for any given $\lambda\neq 0$ there exists a corresponding steady-state solution, ${\sf s}(\lambda)$, in a suitable homogeneous Sobolev  space. Successively, we prove that if there is  $\lambda_\C$ such that the linearization, $\mathscr L$, around ${\sf s}(\lambda_\C)$ is trivial, then there exists a unique analytic family of steady-state solutions in a neighborhood of $\lambda_\C$ where the translational velocity of the sphere is parallel to that of ${\sf s}(\lambda_\C)$; see \theoref{2.2}. In Section 5, we investigate some important spectral properties of  $\mathscr L$ in the Lebesgue space $L^2$. Precisely, we show that the intersection of the spectrum of $\mathscr L$ with the imaginary axis is constituted, at most, by a countable number of eigenvalues of finite multiplicity that can only cluster at 0; see \theoref{5.1}. The main objective of Section 6 is to establish the Fredholm property of the time-periodic linearized operator in a suitable space of functions with zero average over a period. In particular, in \theoref{3.1} we prove that such operator is Fredholm of index 0.  With the help of the results established in the previous sections, in Section 7 we then secure that the original problem is written in an abstract setting where the general theory recalled in Section 3 applies. Therefore, thanks to \theoref{3.1_ar}, in \theoref{7.1} we  give sufficient conditions for the existence of a time-periodic branch in the neighborhood of the steady-state solution ${\sf s}(\lambda_\C)$. As already mentioned, these conditions amount to the request that the operator $\mathscr L$, {\em suitably defined}, has a non-resonating, simple imaginary eigenvalue satisfying the transversality condition. In the final Section 8, we consider the problem of the motion of the sphere in the time-periodic regime. In this regard, we give necessary and sufficient conditions for the occurrence of a horizontal oscillation of the center of mass, in a neigborhood of the ``critical" value $\lambda_\C$; see \theoref{8.1}. 
\section{Formulation of the Problem} A  sphere, $\mathscr S$, of constant density, $\rho_{\mathscr S}$ and radius $R$  freely moves  under the action of gravity in an otherwise quiescent Navier--Stokes liquid, $\mathscr L$,  that fills the entire space outside $\mathscr S$. We assume that $\mathscr S$ is not floating, namely, it has a non-zero buoyancy. This means that, denoting by $\rho_{\mathscr L}$ the  density of the liquid, we take $|\alpha|:=|\rho_{\mathscr S}/\rho_{\mathscr L}-1|>0$. Just to fix the ideas, we shall assume $\alpha>0$ (falling sphere). However, all results continue to hold in the case of positive buoyancy (rising sphere) by simply replacing $\alpha$ with $-\alpha$. 
\par
Let $\mathcal F=\{O,\bfe_1,\bfe_2,\bfe_3\}$ be a frame with the origin at the center of $\mathscr S$ ($\equiv$\, center of mass of $\mathscr S$)  and the axis $\bfe_1$ oriented along the acceleration of gravity $\bfg$. The  dynamics of the coupled system $\mathscr S\cup\mathscr L$ in $\mathcal F$ are then governed by the following set of non-dimensional equations \cite[Section 4]{Gah}
\be\ba{cc}\medskip\left.\ba{ll}\medskip\partial_t\bfv+\lambda\,(\bfv-\bftau)\cdot\nabla\bfv=\Delta\bfv-\nabla p\\
\Div\bfv=0\ea\right\}\ \ \mbox{in $\Omega\times (0,\infty)$}\\ \medskip
\bfv=\bftau+\bfvaromega\times\bfx\ \ \mbox{at $\partial\Omega\times (0,\infty)$}\,,\ \ \Lim{\mbox{\footnotesize $|\bfx|$}\to\infty}\bfv(x,t)=\0\,,\ \ t\in(0,\infty)\,,\\ 
M\dot{\bftau}+\Int{\partial\Omega}{}\mathbb T(\bfv,p)\cdot\bfn=\lambda\,\bfe_1\,,\ \   
\mathcal I\,\dot{\bfvaromega}+\Int{\partial\Omega}{}\bfx\times\mathbb T(\bfv,p)\cdot\bfn=0\, \ \ \mbox{in $(0,\infty)$}\,.
\ea
\eeq{1.1}
Here $\Omega=\real^3\backslash\Omega_0$, with $\Omega_0$ volume occupied by $\mathscr S$.  Furthermore, $\bfv$ and $p+\bfe_1\cdot \bfx$  are (non-dimensional) velocity and  pressure fields of $\mathscr L$, while $\bftau,\bfvaromega$  stand for (non-dimensional) translational and  angular velocities  of $\mathscr S$. Moreover, \be\lambda=\sqrt{\alpha\,g\,R^3}/\nu\,\  (>0)\eeq{Gal} is the (dimensionless) Galilei number, with  $\nu$  kinematic viscosity of $\mathscr L$. Also, $M=4\pi\rho_{\mathscr S}/3\rho_{\mathscr L}$ and $\mathcal I=8\pi\rho_{\mathscr S}/15\rho_{\mathscr L}$ are non-dimensional mass and central moment of inertia of $\mathscr S$, and $\bfn$ is the outer unit normal to $\partial\Omega$.   Finally, 
$$
\mathbb T(\bfv,p)=-p\,\mathds{1}+2\,\mathbb D(\bfv)\,,\ \ \mathbb D(\bfv):=\half\big(\nabla\bfv+(\nabla\bfv)^\top\big)\,,
$$
is the Cauchy tensor with $\mathds{1}$ identity tensor and $\top$ denoting transpose.\par\par
Of particular significance is the subclass of solutions to \eqref{1.1} constituted by those fields $(\bfv_0,p_0,\bftau_0,\bfomega_0)$ that are time independent, namely, they solve the following boundary-value problem
\be\ba{cc}\medskip\left.\ba{ll}\medskip\lambda\,(\bfv_0-\bftau_0)\cdot\nabla\bfv_0=\Delta\bfv_0-\nabla p_0\\
\Div\bfv_0=0\ea\right\}\ \ \mbox{in $\Omega$}\\ \medskip
\bfv_0=\bftau_0+\bfomega_0\times\bfx\ \ \mbox{at $\partial\Omega$}\,,\ \ \Lim{\mbox{\footnotesize $|\bfx|$}\to\infty}\bfv_0(x)=\0\,,\\ 
\Int{\partial\Omega}{}\mathbb T(\bfv_0,p_0)\cdot\bfn=\lambda\,\bfe_1\,,\ \   
\Int{\partial\Omega}{}\bfx\times\mathbb T(\bfv_0,p_0)\cdot\bfn=\0\,.
\ea
\eeq{1.2}\par
Solutions $(\bfv_0,p_0,\bftau_0,\bfomega_0)$ to \eqref{1.2} describe the so called {\em steady free falls} of the sphere in the viscous liquid, and, as explained in the introductory section, their behavior depends on the parameter $\lambda$. 
\par 
In mathematical terms, the time-periodic bifurcation problem can be formulated as follows. Let $\lambda_\C>0$,  let $U=U(\lambda_{\sf c})$ be a neighborhood of $\lambda_{\sf c}$, and let ${\sf s}(\lambda):=(\bfv_0,p_0,\bftau_0,\bfomega_0)(\lambda)$, $\lambda\in U(\lambda_{\sf c})$, be a sufficiently smooth family of solutions to \eqref{1.2}. The objective is then to prove the existence of time-periodic solutions to \eqref{1.1} ``around" {\sf s}$(\lambda_{\sf c})$. Since the period $T:=2\pi/\zeta$ of such solutions is  unknown, it is customary to scale the time by introducing the new variable $\s=\zeta\,t$. Therefore, writing
$$
\bfv(x,t)=\bfu(x,t)+\bfv_0(x)\,,\ \ 
p(x,t)={\sf p}(x,t)+p_0(x)\,,\ \ \bftau(t)=\bfgamma(t)+\bftau_0\,,\ \ \bfvaromega(t)=\bfomega(t)+\bfomega_0
$$
our bifurcation problem means that  we must find a $2\pi$-periodic solution-branch $(\bfu,{\sf p},\bfgamma,\bfomega)(\lambda)$, $\lambda\in U(\lambda_{\sf c})$, to the following set of equations
\be\ba{cc}\medskip\left.\ba{ll}\medskip\zeta\,\partial_\s\bfu-\lambda\bftau_0\cdot\nabla\bfu+\lambda(\bfv_0\cdot\nabla\bfu+(\bfu-\bfgamma)\cdot\nabla\bfv_0)+\lambda(\bfu-\bfgamma)\cdot\nabla\bfu=\Delta\bfu-\nabla{\sf p}\\
\Div\bfu=0\ea\right\}\,  \mbox{in $\Omega\times \real$}\,,\\ \medskip
\bfu=\bfgamma+\bfomega\times\bfx\ \ \mbox{at $\partial\Omega\times \real$}\,,\ \ \Lim{\mbox{\footnotesize $|\bfx|$}\to\infty}\bfu(x,t)=\0\,,\ \ t\in\real\,,\\ 
M\dot{\bfgamma}+\Int{\partial\Omega}{}\mathbb T(\bfu,{\sf p})\cdot\bfn=\0\,,\ \   
\mathcal I\,\dot{\bfomega}+\Int{\partial\Omega}{}\bfx\times\mathbb T(\bfu,{\sf p})\cdot\bfn=\0\, \ \ \mbox{in $\real$}\,.
\ea
\eeq{1.3}
\par
Our strategy to solve this problem  consists in rewriting \eqref{1.3}  as operator equations in suitable Banach spaces, with the involved operators  satisfying a certain number of fundamental properties. Once this goal is accomplished, we will be able to employ the general theory introduced in \cite{GaBif} and recalled in Section \ref{sec:bif}, and derive sufficient conditions for the existence of a time-periodic bifurcating branch.\par
Before performing our study, we recall the main notation used in the paper. With the origin at  the center of $\Omega_0$, we set $B_R:=\{|x|<R\}$, and, for $R>1={\rm diam}\,\Omega_0$, $\Omega_R:=\Omega\cap B_R$,  $\Omega^R:=\Omega\cap\{|x|>R\}$.  As customary, for a domain $A\subseteq\real^3$, $L^q=L^q(A)$ is the Lebesgue space with norm $\|\cdot\|_{q,A}$, and  $W^{m,q}=W^{m,q}(A)$ denotes Sobolev space, $m\in\nat$, $q\in[1,\infty]$, with norm $\|\cdot\|_{m,q,A}$. If $u\in L^q(A)$, $v\in L^{q'}(A)$, $q'=q/(q-1)$, we set $(u,v)_A=\int_Au\,v$. Furthermore, $D^{m,q}=D^{m,q}(A)$ are homogeneous Sobolev spaces with semi-norm $|u|_{m,q,A}:=\sum_{|l|=m}\|D^lu\|_q$. In the above  notation, the subscript $"A"$ will be generally omitted, unless confusion arises. A function $u:A\times \real\mapsto \real^3$ is 
{\em $2\pi$-periodic},  if $u(\cdot,\s+2\pi)=u(\cdot\,\s)$, for a.a. $\s\in \real$,
 and we set
$
{\bar u}:=\frac{1}{2\pi}\int_{0}^{2\pi}u(t){\rm d}t\,.
$
If $B$ is a semi-normed real Banach space with semi-norm $\|\cdot\|_B$, $r=[1,\infty]$,  we denote by  $L^r(0,2\pi;B)$  the class of functions
$u:(0,2\pi)\rightarrow B$ such that 
$$
\|u\|_{L^r(B)}\equiv\left\{\ba{ll}\smallskip\big( \Int{0}{2\pi}\|u(t)\|_B^r {\rm d}t\big)^{\frac 1r}<\infty, \ \ \mbox{if 
$r\in [1,\infty)\,;$}\\   
\essup{t\in[0,2\pi]}\,\|u(t)\|_B <\infty, \ \ \mbox{if $r=\infty.$}
\ea\right.
$$
Furthermore, we define
$$
W^{1,2}(0,2\pi;B)=\Big\{u\in L^{2}(0,2\pi;B): \textcolor{black}{\partial_tu\in L^{2}(0,2\pi;B)\,,}\Big\}\,.
$$
Unless otherwise stated, we shall  write $L^r(B)$ for $L^r(0,2\pi;B)$, etc. By $B_{\mathbb C}:=B+{\rm i}\, B$ we denote the complexification of $B$.
If $M$ is a map between two spaces,  ${\sf D}\,[M]$, ${\sf N}\,[M]$, ${\sf R}\,[M]$ and ${\sf P}\,[M]$ will indicate its domain, null space,  range, and resolvent set, respectively.
Finally, by  $c$, $c_0$, $c_1$, etc.,  we denote positive constants, whose particular value is unessential to the context. When we wish to emphasize
the dependence of $c$ on some parameter $\mu$, we shall write  $c(\mu)$ or $c_\mu$.
\setcounter{equation}{0}
\section{An Abstract Bifurcation Theorem}\label{sec:bif}
Objective of this section is to recall a  time-periodic bifurcation theorem for a general class of  operator equations proved in \cite{GaBif}.  Before stating the result, however, we first would like to make some comments that will also provide the motivation of this approach.

Many evolution problems in mathematical physics can be formally written in the form
\be
u_t+L(u)=N(u,\mu)\,,
\eeq{Ar.1}
where $L$ is a linear differential operator (with appropriate {\em homogeneous} boundary conditions), and $N$ is a nonlinear operator depending on the parameter $\mu\in\real$, such that $N(0,\mu)=0$ for all admissible values of $\mu$. Then, roughly speaking, time-periodic bifurcation for  \eqref{Ar.1} amounts to show the existence a family of non-trivial time-periodic solutions $u=u(\mu;t)$ of (unknown) period $T=T(\mu)$ ($T$-{\em periodic} solutions) in a neighborhood of $\mu=0$, and such that $u(\mu;\cdot)\to 0$ as $\mu\to 0$. Setting $\s:=2\pi\,t/T\equiv \zeta\, t$, \eqref{3.1} becomes
\be 
\zeta\,u_\s+L(u)=N(u,\mu)
\eeq{Ar.2}
and the problem reduces to find a family of $2\pi$-periodic solutions to \eqref{Ar.2} with the above properties. We now write $u=\bar{u}+(u-\bar{u}):=v+w$ and observe that \eqref{Ar.2} is formally equivalent to the following two equations
\be\ba{ll}\medskip 
L(v)=\bar{N(v+w,\mu)}:=N_1(v,w,\mu)\,,\\ \zeta\,w_\s+L(w)=N(v+w,\mu)-\bar{N(v+w,\mu)}:=N_2(v,w,\mu)\,.\ea
\eeq{Ar.3}
At this point, the crucial issue is that in many applications --typically when the physical system evolves in an {\em unbounded spatial region}-- the ``steady-state component" $v$ lives in function spaces with quite less ``regularity''\footnote{Here `regularity' is meant in the sense of behavior at large spatial distances.} than the space where the ``oscillatory" component $w$ does. For this reason, it is much more appropriate to study the two equations in \eqref{Ar.3} in two {\em different} function classes. As a consequence, even though {\em formally}  being the same as differential operators, the operator $L$ in \eqref{Ar.3}$_1$ acts on and ranges into spaces  different than those the operator $L$ in \eqref{Ar.3}$_2$ does. With this in mind, \eqref{Ar.3} becomes
$$
L_1(v)=N_1(v,w,\mu)\,;\ \ \zeta\,w_\s+L_2(w)=N_2(v,w,\mu)\,.
$$   
\par
The general abstract theory that we are about to describe stems exactly from the above considerations. 
\par  
Let $\mathcal X, \mathcal Y$,  be (real) Banach spaces with norms $\|\cdot\|_{\mathcal X}$, $\|\cdot\|_{\mathcal Y}$, respectively, and let  $\mathcal Z$ be a (real) Hilbert space with norm $\|\cdot\|_{\mathcal Z}$ and corresponding scalar product $\langle\cdot,\cdot\rangle$. Moreover,  denote by 
$$
L_1:\mathcal X\mapsto \mathcal Y\,,
$$
a bounded linear operator, and by
$$
L_2:{\sf D}\,[L_2]\subset \mathcal Z\mapsto \mathcal Z\,,
$$
a densely defined, closed linear operator, with a non-empty resolvent set ${\sf P}(L_2)$. For a fixed (once and for all) $\theta\in {\sf P}(L_2)$ we denote by $\mathcal W$ the linear subspace of $\mathcal Z$ closed under the norm $\|w\|_{\mathcal W}:=\|(L_2+\theta\,I)w\|_{\mathcal Z}$, where $I$ stands for the identity operator. We also define the following spaces
$$\ba{rl}\medskip
\mathcal Z_{2\pi,0}&\!\!\!:=\big\{w\in L^2(0,2\pi;\mathcal Z): \mbox{$2\pi$-periodic with} \ \bar{w}=0\big\} 
\\\medskip\calw_{2\pi,0}&\!\!\!:=\big\{w\in L^2(0,2\pi;\calz)\,, \  w_t\in L^{2}(0,2\pi;\mathcal Z): \mbox{$2\pi$-periodic with} \ \bar{w}=0\big\}\,.
\ea$$ 
Next, let
$$
N: \calx\times \calw_{2\pi,0}\times \real\mapsto \caly\oplus \calh_{2\pi,0}   
$$
be a (nonlinear) map satisfying the following properties:
\be\ba{rl}\medskip
N_1&\!\!\!: (v,w,\mu)\in\calx\times \calw_{2\pi,0}\times \real\mapsto \bar{N(v,w,\mu)}\in \caly
\\
N_2&\!\!\!:=N-N_1:\calx\times \calw_{2\pi,0}\times \real\mapsto \calz_{2\pi,0}\,.
\ea
\eeq{Ar.4} 
The bifurcation problem can be then rigorously formulated as follows.\smallskip\par 
Bifurcation Problem: {\em Find a neighborhood of the origin $U(0,0,0)\subset \calx\times \calw_{2\pi,0}\times \real$ such that the equations
\be
L_1(v)=N_1(v,w,\mu)\,,\ \mbox{in $\caly$}\,;\ \ \zeta\, w_\s +L_2(w)=N_2(v,w,\mu)\,,\ \mbox{in $\calz_{2\pi,0}$}\,,
\eeq{Ar.5}
possess there a family of non-trivial $2\pi$-periodic solutions $(v(\mu),w(\mu;\tau))$ for some $\omega=\omega(\mu)>0$,  such that $(v(\mu),w(\mu;\cdot))\to 0$ in $\calx\times\calw_{2\pi,0}$ as $\mu\to0$.}\smallskip\par
Whenever the Bifurcation Problem admits a positive answer, we say that  $(u=0,\mu=0)$ is a {\em bifurcation point}. Moreover, the bifurcation is called {\em supercritical} [resp. {\em subcritical}] if the family of solutions $(v(\mu),w(\mu;\tau))$ exists only for $\mu>0$ [resp. $\mu<0$]. 
\medskip\par
In order to provide sufficient conditions for the resolution of the above problem, we need the following  assumptions (H1)--(H4) on the involved operators.
\begin{itemize}
  \item[(H1)] $L_1$ is a homeomorphism\,;
  \item[(H2)] There exists $\nu_0:={\rm i}\,\zeta_0$, $\zeta_0>0$ such that   $L_2-\nu_0I$ is Fredholm of index 0,  and   ${\rm dim}\,{\sf N}_{\mathbb C}[L_2-\nu_0I]=1$ with ${\sf N}_{\mathbb C}[L_2-\nu_0I]\cap {\sf R}_{\mathbb C}[L_2-\nu_0I]=\{0\}$. Namely, $\nu_0$ is a simple  eigenvalue of $L_2$. Moreover,  $k\,\nu_0\in {\sf P}(L_2)$, for all $k\in\nat\backslash\{0,1\}$. 
  \item[(H3)] The operator
$$
\mathscr Q:w\in \calw_{2\pi,0}\mapsto \zeta_0\,w_\s+L_2(w)\in \calz_{2\pi,0}\,,
$$  
is Fredholm of index 0\,;
\item[(H4)] The nonlinear operators $N_1,N_2$
are analytic in the neighborhood $U_1(0,0,0)\subset \calx\times \calw_{2\pi,0}\times \real$, namely, there exists $\delta>0$ such that for all $(v,w,\mu)$ with $\|v\|_{\calx}+\|w\|_{\calw_{2\pi,0}}+|\mu|<\delta$, the Taylor series 
$$\ba{ll}\medskip
N_1(v,w,\mu)=\Sum{k,l,m=0}{\infty}R_{klm}v^kw^l\mu^m\,,\\
N_2(v,w,\mu)=\Sum{k,l,m=0}{\infty}S_{klm}v^kw^l\mu^m\,,
\ea
$$
are absolutely convergent in $\caly$ and $\calz_{2\pi,0}$, respectively, for all $(v,w,\mu)\in U_1$. Moreover, we assume that the multi-linear operators $R_{klm}$ and $S_{klm}$ satisfy $R_{klm}=S_{klm}=0$ whenever $k+l+m\le1$, and $R_{011}=R_{00m}=S_{00m}=0$, all $m\ge2$.
\end{itemize}
\par
Before stating the bifurcation result, we need to recall a relevant consequence of the above assumptions. Let
$$
L_2(\mu):=L_2+\mu\,S_{011}\,,
$$
and observe that, by (H2), $\nu_0$ is a simple eigenvalue of $L_2(0)\equiv L_2$. Therefore, denoting by $\nu(\mu)$ the eigenvalues of $L_2(\mu)$, it follows (e.g.  \cite[Proposition 79.15 and Corollary 79.16]{Z1}) that in a neighborhood of $\mu=0$ the map $\mu\mapsto\nu(\mu)$ is well defined and of class $C^\infty$.
\renewcommand{\theequation}{\arabic{section}.\arabic{equation}}
\smallskip\par\setcounter{equation}{5}
We may then state the following bifurcation result, whose proof is given in \cite[Theorem 3.1]{GaBif}.
\Bt Suppose  {\rm (H1)--(H4)} hold and, in addition,
$$
\Re[\nu'(0)]\neq 0\,,
$$
namely, the  eigenvalue $\nu(\mu)$ crosses the imaginary axis with ``non-zero speed." Moreover, let $v_0$ be a normalized eigenvector of $L_2$ corresponding to the eigenvalue $\nu_0$, and set
$
v_1:=\Re[v_0\,{\rm e}^{-{\rm i}\,\s}].$    
Then, the following properties are valid. \smallskip\\
{\rm (a)} {\rm Existence.} There are analytic families
\be
\big(v(\varepsilon),w(\varepsilon),\zeta(\varepsilon),\mu(\varepsilon)\big)\in \calx\times \calw_{2\pi,0}\times \real_+\times\real
\eeq{fam}
satisfying \eqref{Ar.5}, for all real $\varepsilon$ in a neighborhood $\mathcal I(0)$ of 0, and such that
\be
\big(v(\varepsilon),w(\varepsilon)-\varepsilon\,v_1,\zeta(\varepsilon),\mu(\varepsilon)\big)\to (0,0,\zeta_0,0)\ \ \mbox{as $\varepsilon\to 0$}\,.
\eeq{Ar.10}
\par\noindent
{\rm (a)} {\rm Uniqueness.}
There is a neighborhood  $$U(0,0,\zeta_0,0)\subset \calx\times \calw_{2\pi,0}\times \real_+\times \real$$ such that every (nontrivial) $2\pi$-periodic solution to \eqref{Ar.5},  lying in $U$ must coincide, up to a phase shift, with a member of the family \eqref{fam}.
\smallskip\par\noindent
{\rm (a)} {\rm Parity.}  The functions $\zeta(\varepsilon)$ and $\mu(\varepsilon)$ are even:
$$
\zeta(\varepsilon)=\zeta(-\varepsilon)\,,\ \ \mu(\varepsilon)=\mu(-\varepsilon)\,,\ \ \mbox{for all $\varepsilon\in\cali(0)$\,.} 
$$
Consequently, the bifurcation due to these solutions is either subcritical or supercritical, a two-sided bifurcation being excluded.\footnote{Unless $\mu\equiv 0$.}
\ET{3.1_ar}

 \setcounter{equation}{0}
\section{Relevant Function Spaces and related Properties}
In this section we will introduce  certain function classes along with some of their most important  properties. 
\par
Let $\bfpzc R$ be the class of velocity fields in a rigid motion:
$$
\bfpzc R=\{\hat{\bfu}\in C^\infty(\real^3):\ \hat{\bfu}=\hat{\bfu}_1+\hat{\bfu}_2\times\bfx\,, \  \mbox{for some $\hat{\bfu}_1,\hat{\bfu}_2\in\real^3$}\}\,,
$$
and set
$$\ba{ll}\medskip
\calk=\mathcal K(\real^3)=\big\{\bfphi\in C_0^\infty(\bar{\real^3}):\ \mbox{in $\real^3$}\,;
\bfphi(x)=\hat{\bfphi}, \ \mbox{some $\hat{\bfphi}\in\bfpzc R$}\, \ \mbox{in a neighborhood of $\Omega_0$}\big\}
\\
\calc=\mathcal C(\real^3)=\{\bfphi\in\calk(\real^3):\ \Div\bfphi=0\ \mbox{in $\real^3$}\}\,.
\ea
$$
We shall call $\hat{\bfphi}_1,\hat{\bfphi}_2$ the {\em characteristic vectors} of the function $\bfphi$.
In $\calk$ we introduce the scalar product
\be
M\hat{\bfphi}_1\cdot\hat{\bfpsi}_1+\cali \hat{\bfphi}_2\cdot \hat{\bfpsi}_2+(\bfphi,\bfpsi)_\Omega\,,\ \ \bfphi,\bfpsi\in\calk\,,
\eeq{0.0}
and define the following spaces\Footnote{Even though $\bfx\times\bfn=\0$ if $x\in\partial\Omega$, for notational convenience we will keep the term $\hat{\bfh}_2$.}
\be\ba{ll}\medskip
\mathcal L^2=\call^2(\real^3):= \,\big\{\mbox{completion of $\calk(\real^3)$ in the norm induced by \eqref{0.0}}\big\}\,,\\ \medskip
\mathcal H=\calh(\real^3):=\,\big\{\mbox{completion of $\calc(\real^3)$ in the norm induced by \eqref{0.0}}\big\}\\ \medskip
\calg=\calg(\real^3):=\big\{\bfh\in \call^2(\real^3): \, \mbox{there is $p\in D^{1,2}(\Omega)$ such that $\bfh=\nabla p$ in $\Omega$}\,,\\ \hspace*{3cm} \mbox{and}\ \bfh=-M^{-1}\int_{\partial\Omega}p\,\bfn   -\cali^{-1}\left(\int_{\partial\Omega}p\bfy\times\,\bfn\right)\times\bfx\ \mbox{in}\ \Omega_0\,\big\}\,.
\ea
\eeq{spazi}
In \cite[Theorem 3.1 and Lemma 3.2]{ALS} the following characterization of the spaces $\call^2$ and $\calh$ is proved.
\Bl Let $\Omega$ be Lipschitz. Then
$$\ba{ll}\medskip
\call^2(\real^3)=\{\bfu\in L^2(\real^3): \ \bfu=\hat{\bfu}\ \mbox{in}\ \Omega_0, \ \mbox{for some $\hat{\bfu}\in \bfpzc R$}\}\\
\calh(\real^3)=\{\bfu\in \call^2(\real^3): \ \Div\bfu=0\,\}\,.
\ea
$$
\EL{aria}
Furthermore, by an argument entirely analogous to that employed in \cite[Theorem 3.2]{ALS} one shows:
\Bl The following orthogonal decomposition holds
$$
\call^2(\real^3)=\calh(\real^3)\oplus\calg(\real^3)\,.
$$
\EL{0.1}
We next introduce the space
$$
\calv=\calv(\real^3)\equiv \,\big\{\mbox{completion of $\calc(\real^3)$ in the norm $\|\mathbb D(\cdot)\|_2$}\big\}\,.
$$
The basic properties of the space $\calv$ are collected in the next lemma, whose proof is given in \cite[Lemmas 9--11]{Gah}.
\Bl  $\calv$ is a Hilbert space endowed with the scalar product
\be
\big[\bfu,\bfw\big]=\int_\Omega\mathbb D(\bfu):\mathbb D(\bfw)\,,\ \ \bfu,\bfw\in\calv\,.
\eeq{1.6}
Furthermore, the following characterization holds
\be
\calv=\big\{\bfu\in L^6(\real^3)\cap D^{1,2}(\real^3)\,;\ \Div\bfu=0\ \mbox{in\, $\real^3$}\,;
\bfu(y)=\hat{\bfu}\,,\ y\in\Omega_0,  \ \mbox{some $\hat{\bfu}\in\bfpzc R$} \big\}\,.
\eeq{1.7}
Also, for each $\bfu\in\calv$, we have
\be
\|\nabla\bfu\|_2=\sqrt{2}\|\mathbb D(\bfu)\|_2\,,
\eeq{1.8}
and
\be
\|\bfu\|_6\le \kappa_0\,\|\mathbb D(\bfu)\|_2\,,\ \ \bfu\in\calv\,,
\eeq{1.9}
for some numerical constant $\kappa_0>0$.\footnote{Recall that, in our non-dimensionalization, the sphere has radius 1.} Finally, 
there is another positive numerical constant $\kappa_1$ such that
\be
|\hat{\bfu}_1|+|\hat{\bfu}_2|\le \kappa_1\,\|\mathbb D(\bfu)\|_2\,.
\eeq{1.10}
\EL{1.1}
\par
Let $\calv^{-1}=\calv^{-1}(\real^3)$ be the dual space of $\calv$, and let  $\langle\cdot,\cdot\rangle$ and $\|\cdot\|_{-1}$ be the corresponding duality pair and associated norm, respectively. Denote by $\bfe\in\real^3$ a given unit vector and consider the space
$$
\mathscr X_\sub{\bfe}(\real^3):=\big\{\bfu\in\calv(\real^3): \ \bfe\cdot\nabla\bfu\in\calv^{-1}\big\}
$$ 
where $\bfe\cdot\nabla\bfu\in\calv^{-1}$ means that there is $C=C(\bfu)>0$ such that
\be
|(\bfe\cdot\nabla\bfu,\bfphi)|\le C\,\|\mathbb D(\bfphi)\|_2\,,\ \ \mbox{for all $\bfphi\in\calc(\real^3)$}\,.
\eeq{1.11}
Actually, from \eqref{1.11},  the density of $\calc$ in $\calv$, and the Hahn--Banach theorem it follows that 
$\bfe\cdot\nabla\bfu$ can be uniquely extended to a bounded linear functional on the whole of $\calv$, with
$$
\|\bfe\cdot\nabla\bfu\|_{-1}:=\sup_{\mbox{\footnotesize $\bfphi\in\calc; \|\mathbb D(\bfphi)\|_2=1$}}|(\bfe\cdot\nabla\bfu,\bfphi)|\,.  
$$
Obviously, the functional 
$$
\|\bfu\|_{{\mathscr X}_\sub{\bfe}}:=\|\mathbb D(\bfu)\|_2+\|\partial_1\bfu\|_{-1}\,,
$$
defines a norm in $\mathscr X(\real^3)$. 

In  the following lemma we collect the relevant properties of the space $\mathscr X(\real^3)$. Their proofs are entirely analogous to \cite[Proposition 65]{Gace} and \cite[Lemma 2.1]{GaJMP}, and therefore will be omitted.
\Bl The space
$\mathscr X_\sub{\bfe}(\real^3)$ endowed with the norm $\|\cdot\|_{{\mathscr X}_\sub{\bfe}}$ is a reflexive, separable Banach space, dense in $\calv(\real^3)$.
Moreover, $\mathscr X_\sub{\bfe}(\real^3)$ is continuously embedded in $L^4(\real^3)$, and there is $c>0$ such that
$$
\|\bfu\|_4\le c\,\big(\|\bfe\cdot\nabla\bfu\|_{-1}^{\frac14}\|\mathbb D(\bfu)\|_2^{\frac34}+\|\mathbb D(\bfu)\|_2\big)\,.
$$
Finally, we have
\be
\langle\bfe\cdot\nabla\bfu,\bfu\rangle=0\,,\ \ \mbox{for all $\bfu\in \mathscr X_\sub{\bfe}(\real^3)$.}
\eeq{1.12}
\EL{1.2}
\par
We next introduce the following spaces of time-periodic functions: 
$$\ba{ll}\medskip
L^{2}_\sharp=\{\bfF\in L^{2}(0,2\pi); \ \mbox{$\bfF$ is $2\pi$-periodic with }\ \bar{\bfF}=0\}\\\medskip
W^{1,2}_\sharp=\{\bfchi\in L^{2}_\sharp(0,2\pi); \ \dot{\bfchi}\in L^2(0,2\pi)\}\\
\medskip
\mathcal L_\sharp^{2}:=\{\bff\in L^{2}(L^2); \ \mbox{$\bff$ is $2\pi$-periodic, with $\bar{\bff}=\0$}\}
\\ 
\mathcal W_\sharp^{2}:=\{\bfu\in W^{1,2}(L^2)\cap L^2(W^{2,2}\cap\calv); \ \mbox{$\bfu$ is $2\pi$-periodic, with $\bar{\bfu}=\0$}\}
\ea
$$
endowed with  natural norms
$$\ba{ll}\medskip
\|\bfchi\|_{L^{2}_\sharp}:=\|\bfchi\|_{L^2}\,;\ \ \
\|\bfchi\|_{{W}_\sharp^{1,2}}:=\|\bfchi\|_{{W}^{1,2}}\equiv \|\dot{\bfchi}\|_{L^{2}}+\|\bfchi\|_{L^2}\,,\\
\|\bfu\|_{\mathcal L^{2}_\sharp}:=\|\bfu\|_{L^2(L^2)}\,;\ \ \
\|\bfu\|_{{\mathcal W}_\sharp^{2}}:=\|\bfu\|_{W^{1,2}(L^2)}+\|\bfu\|_{L^2(W^{2,2})}\,.
\ea
$$
Finally, we set
$$
\calp^{1,2}:= L^2(D^{1,2})\,.
$$

\setcounter{equation}{0}
\section{On the Existence of a Steady-State Branch}
The goal of this section is two-fold. On the one hand, to establish existence of solutions to \eqref{1.2} in suitable function classes for all values of $\lambda>0$ and, on the other hand, to furnish sufficient conditions for the existence of a smooth solution branch around the ``critical" value $\lambda_\C$. 
\par
We begin to give the definition of weak solution.
\Bd A field $\bfv:\real^3\mapsto \real^3$ is a {\em weak solution} to \eqref{1.2} if the following conditions hold.\begin{itemize}\item[{\rm (i)}] $\bfv\in \calv(\real^3)$, and $\bftau\equiv\hat{\bfv}_1\neq\0$\,;
\item[{\rm (ii)}] $\bfv$ satisfies the relation:\footnote{For simplicity, we omit the subscript ``${}_\sub{0}"$.}
\be
2[\bfv, \bfphi]-\lambda\,\hat{\bfphi}_1\cdot\bfe_1=\lambda\big((\bftau-\bfv)\cdot\nabla\bfv,\bfphi\big)\,,\ \ \mbox{for all $\bfphi\in\calc$}\,.
\eeq{2.1}
\end{itemize} 
\EDD{2.1}
\Br It is easy to show that,
 if a weak solution $\bfv$ is sufficiently smooth in $\Omega$, then there exists a likewise smooth pressure field $p:\Omega\mapsto\real$ such that the quadruple $(\bfv,p,\bftau\equiv\hat{\bfv}_1,\bfomega\equiv\hat{\bfv}_2)$
is a smooth solution to \eqref{1.2}. To this end, we begin to obseve that, clearly, the validity of \eqref{1.2}$_{2,3}$ comes from the fact that $\bfv\in\calv(\real^3)$. If we integrate by parts the first term on the left-hand side of \eqref{2.1}, we get
\be
2[\bfv, \bfphi]=-\Int{\Omega}{}\Delta\bfv\cdot\bfphi+\hat{\bfphi}_1\cdot\int_{\partial\Omega}2\,\mathbb D(\bfv)\cdot\bfn+\hat{\bfphi}_2\cdot\int_{\partial\Omega}\bfx\times(2\,\mathbb D(\bfv)\cdot\bfn)\,.
\eeq{ar}
Thus, by taking, in particular, $\bfphi\in \cald(\Omega):=\{\bfphi\in C_0^\infty(\Omega): \Div\bfphi=0\}$, from \eqref{2.1} and \eqref{ar} we get 
$$
\int_{\Omega}(-\Delta\bfv-\lambda(\bftau-\bfv)\cdot\nabla\bfv)\cdot\bfphi=0\,,\ \ \mbox{for all $\bfphi\in\cald(\Omega)$\,,}
$$
which, by well-known results, implies  the existence of a (smooth) scalar field $p$ for which \eqref{1.2}$_{1}$ holds. Now, using \eqref{ar} and \eqref{1.2}$_1$ we deduce
$$
\hat{\bfphi}_1\cdot\big(-\lambda\,\bfe_1+\int_{\partial\Omega}2\,\mathbb D(\bfv)\cdot\bfn\big)+\hat{\bfphi}_2\cdot\int_{\partial\Omega}\bfx\times(2\,\mathbb D(\bfv)\cdot\bfn)-\int_\Omega\nabla p\cdot\bfphi=\lambda \int_{\Omega_0}(\bftau-\bfv)\cdot\nabla\bfv\cdot\hat{\bfphi}\,,
$$
which, after an integration by parts leads to
\be
\hat{\bfphi}_1\cdot\big(-\lambda\,\bfe_1+\int_{\partial\Omega}\mathbb T(\bfv,p)\cdot\bfn\big)+\hat{\bfphi}_2\cdot\int_{\partial\Omega}\bfx\times\mathbb T(\bfv,p)\cdot\bfn=\lambda \int_{\Omega_0}(\bftau-\bfv)\cdot\nabla\bfv\cdot\hat{\bfphi}\,.
\eeq{ar1}
Recalling that, in $\Omega_0$, $\bfv,\bfphi\in \bfpzc R$, we get ($\bftau\equiv\hat{\bfv}_1$)
$$\ba{rl}\medskip
\Int{\Omega_0}{}(\bftau-\bfv)\cdot\nabla\bfv\cdot\hat{\bfphi}&\!\!\!\!=\Int{\Omega_0}{}(\hat{\bfv}_2\times\bfx)\cdot\nabla(\hat{\bfv}_2\times\bfx)\cdot\hat{\bfphi}\\
&\!\!\!\!=\Int{\partial\Omega}{}\hat{\bfv}_2\times\bfx\cdot\bfn (\hat{\bfv}_2\times\bfx\cdot\hat{\bfphi})-\Int{\Omega}{}(\hat{\bfv}_2\times\bfx)\cdot\mathbb D(\hat{\bfphi})\cdot(\hat{\bfv}_2\times\bfx)=0
\ea
$$
because $\mathbb D(\hat{\bfphi})=0$, and $\bfx\times\bfn=\0$ for all $x\in\partial\Omega$. From the latter, \eqref{ar1}, and the arbitrariness of $\hat{\bfphi}$ we then conclude that a smooth weak solution $\bfv$ and the associated pressure field $p$ satisfy also \eqref{1.2}$_{5,6}$. Finally, as a consequence of the next lemma, we shall see that $\bfv$ obeys also the asymptotic condition \eqref{1.2}$_4$.
\ER{2.1}
\par
The next lemma establishes some important properties of weak solutions.
\Bl Let $\bfv$ be a weak solution to \eqref{1.2}. Then, the following properties hold.
\begin{itemize}
\item[{\rm (a)}] $\bfv\in C^\infty(\bar{\Omega}_R)\cap C^\infty(\Omega)$,  all $R>1$,  and there is a pressure field $p\in C^\infty(\bar{\Omega}_R)\cap C^\infty(\Omega)$,  all $R>1$, such that the quadruple $(\bfv,p,\bftau\equiv\hat{\bfv}_1,\bfomega\equiv \hat{\bfv}_2)$  satisfies \eqref{1.2}\,;
\item[{\rm (b)}] $\bfv\in L^q(\real^3)\cap D^{1,\frac{2q}3}(\real^3)$\,,\ $p\in L^{\frac {3q}4}(\Omega)$\,, for all $q\in (2,\infty]$\,; 
\item[{\rm (c)}]  $\bfv\in \mathscr X_\sub{\bfe}(\real^3)$\,,\ $\bfe:=\bftau/|\bftau|$\,;  
\item[{\rm (d)}] $\bfv\cdot\nabla\bfv\in\calv^{-1}(\real^3)$ and relation \eqref{2.1} is equivalent to  the following one:
\be
2[\bfv, \bfphi]-\lambda\,\hat{\bfphi}_1\cdot\bfe_1=\lambda\big\langle(\bftau-\bfv)\cdot\nabla\bfv,\bfphi\big\rangle\,,\ \ \mbox{for all $\bfphi\in\calv(\real^3)$}\,.
\eeq{2.6}
\item[{\rm (e)}] $\bfv$ obeys the energy equality:
\be
2\|\mathbb D(\bfv)\|_2^2=\lambda\,\bftau\cdot\bfe_1\,.
\eeq{2.2}
\end{itemize}
\EL{2.1}
{\em Proof.} Taking into account that $\Omega$ is of class $C^\infty$, the first statement, restricted to \eqref{1.2}$_{1,2,3,5,6}$, follows from the regularity results in  \cite[Theorem IX.5.1]{Gab} applied to \eqref{2.1}, along with \remref{2.1}. Next, since $\bftau\neq\0$ and $\bfv\in\calv(\real^3)$, the results in  \cite[Theorem X.6.4]{Gab} ensure that $\bfv$, and $p$ possess  the summability properties stated in (b). In particular, by \cite[Theorem II.9.1]{Gab}, the latter imply also \eqref{1.2}$_4$.
By H\"older  inequality and \eqref{1.8}, \eqref{1.9} we show, for all $\bfphi\in\calc$,
\be
|(\bfv\cdot\nabla\bfv,\bfphi)|\le \|\bfv\|_3\|\nabla\bfv\|_2\|\bfphi\|_6\le \sqrt{2}\kappa_0\|\bfv\|_3\|\mathbb D(\bfv)\|_2\|\mathbb D(\bfphi)\|_2\,,
\eeq{2.4}
which delivers in particular, by (b), \be\bfv\cdot\nabla\bfv\in\calv^{-1}(\real^3).\eeq{2.3_0} 
Furthermore, from \eqref{2.1}, \eqref{1.10}  and \eqref{2.4} we infer, for all $\bfphi\in\calc$, 
$$\ba{rl}\medskip
|(\bftau\cdot\nabla\bfv,\bfphi)|\le &\!\!\!\!\|\mathbb D(\bfv)\|_2\|\mathbb D(\bfphi)\|_2+\lambda\,\big[|\hat{\bfphi}_1|+\|\bfv\|_3\|\nabla\bfv\|_2\|\bfphi\|_6\big]\\
\le &\!\!\!\!\big[\|\mathbb D(\bfv)\|_2+\lambda\,(\kappa_1+\sqrt{2}\kappa_0\|\bfv\|_3\|\mathbb D(\bfv)\|_2)\big]\,\|\mathbb D(\bfphi)\|_2
\,,
\ea$$
from which we conclude $\bfv\in\mathscr X_\sub{\bfe}(\real^3)$.
In view of the latter,  \eqref{2.3_0}, and the density of $\calc$ in $\calv(\real^3)$, we may thus deduce that  \eqref{2.1} leads to \eqref{2.6}.
We now replace $\bfv$ for $\bfphi$ in \eqref{2.6} to obtain
\be
2\|\mathbb D(\bfv)\|_2^2-\lambda\,\bftau\cdot\bfe_1=\lambda\big(\langle\bftau\cdot\nabla\bfv,\bfv\rangle-\langle\bfv\cdot\nabla\bfv,\bfv\rangle\big)=-\lambda\,\langle\bfv\cdot\nabla\bfv,\bfv\rangle\,,
\eeq{A1}
where, in the last step, we have taken into account \eqref{1.12}. Now, after integrating by parts,
we deduce 
\be
(\bfv\cdot\nabla\bfv,\bfphi)=-(\bfv\cdot\nabla\bfphi,\bfv)\,,\ \ \mbox{for all $\bfphi\in\calc$.}
\eeq{int}
Since $\bfv\in L^4(\real^3)$ by the 
H\"older inequality it is readily seen that the right-hand side of \eqref{int} defines a continuous functional in $\bfphi\in\calv(\Omega)$. As a result, setting $\bfphi=\bfv$ in \eqref{int} we conclude $(\bfv\cdot\nabla\bfv,\bfv)=\langle\bfv\cdot\nabla\bfv,\bfv\rangle=0$, and the property (e) follows from this and \eqref{A1}.
The proof of the lemma is completed.  
\par\hfill$\square$\par
We are now in a position to show the following existence result.
\Bt
For any given $\lambda>0$, problem \eqref{1.2} has at least one weak solution $\bfv=\bfv(\lambda)$.
\ET{2.1}
{\em Proof.} We shall employ the classical Galerkin method.  To this end, let $\{\bfphi_k\}\subset\calc$ be an ortho-normal basis in $\calv(\real^3)$ and, for notational simplicity, set 
$$
\widehat{(\bfphi_k)}_1=\bfxi_k\,,\ \ \ \widehat{(\bfphi_k)}_2=\bfsigma_k\,.  
$$
Consider the linear combinations
$$ 
\bfv_m=\sum_{\ell=1}^mc_{\ell m}\bfphi_\ell\,,\ \ \bftau_m=\sum_{\ell=1}^mc_{\ell m}\bfxi_\ell\,,\ \ \bfomega_m=\sum_{\ell=1}^mc_{\ell m}\bfsigma_\ell
$$
where the coefficients $c_{\ell m}$ are requested to be solution to the following algebraic set of equations
\be
[\bfv_m, \bfphi_k]-\lambda\,\bfxi_k\cdot\bfe_1=\lambda\big((\bftau_m-\bfv_m)\cdot\nabla\bfv_m,\bfphi_k)\,,\ \ \mbox{$k\in\{1,\ldots m\}$}\,.
\eeq{2.8} 
Multiplying both sides of this relation by $c_{k m}$ and summing over $k\in\{1,\ldots,m\}$ we show
\be
2\|\mathbb D(\bfv_m)\|_2^2-\lambda\bftau_m\cdot\bfe_1=\lambda\,\big((\bftau_m-\bfv_m)\cdot\nabla\bfv_m,\bfv_m\big)\,.
\eeq{2.9}
By a simple integration by parts,  we show
$$
\big((\bftau_m-\bfv_m)\cdot\nabla\bfv_m,\bfv_m\big)=0\,,
$$ 
and so \eqref{2.9} becomes
\be
2\|\mathbb D(\bfv_m)\|_2^2-\lambda\bftau_m\cdot\bfe_1=0\,.
\eeq{2.10}
Now, let $\bfc:=(c_{1m},\ldots,c_{mm})\in\real^m$, and consider the map
$$
\bfP:\bfc\in\real^m\mapsto \bfP(\bfc)\in\real^m\,,
$$
where
$$
[\bfP(\bfc)]_k=[\bfv_m, \bfphi_k]-\lambda\,\bfxi_k\cdot\bfe_1-\lambda\big((\bftau_m-\bfv_m)\cdot\nabla\bfv_m,\bfphi_k)\,,\ \ \mbox{$k\in\{1,\ldots m\}$}\,.
$$
By what we have just shown,
$$
\bfP(\bfc)\cdot\bfc=2\|\mathbb D(\bfv_m)\|_2^2-\lambda\bftau_m\cdot\bfe_1\,,\ \ \mbox{for all $m\in\nat$},
$$
and so, by \eqref{1.10} and the ortho-normal property of the base $\{\bfphi_k\}$,
$$
\bfP(\bfc)\cdot\bfc\ge |\bfc|(|\bfc|-\kappa_1)
$$ 
Thus, by \cite[Lemma IX.3.1]{Gab}, the latter implies that for each $m\in\nat$, the algebraic system \eqref{2.8} has at least one solution $\bfc=\bfc(m)$. Furthermore, from \eqref{2.10} and again \eqref{1.10} we deduce the following estimate for the sequence $\{\bfv_m\}$, uniformly in $m\in\nat$:
\be
2\|\mathbb D(\bfv_m)\|_2\le \lambda\,\kappa_1\,.
\eeq{2.11}
Taking into account that $\calv(\real^3)\subset W^{1,2}(B_R)$ and that $W^{1,2}(B_R)$ is compactly embedded in $L^2(B_R)$, for all $R>0$, from \eqref{2.11} we deduce the existence of $\bfv\in\calv(\real^3)$  such that (possibly along a subsequence)
$$
\ba{ll}\medskip
\bfv_m\to \bfv\,,\ \mbox{weakly in $\calv(\real^3)$,\, strongly in $L^2(B_R)$, all $R>0$}\\
\bftau_m\to \bftau\equiv\hat{\bfv}_1,\ \, \bfomega_m\to \bfomega\equiv\hat{\bfv}_2\,\ \mbox{in $\real^3$}\,.
\ea
$$
Employing these convergence properties, we may pass to the limit $m\to\infty$ in \eqref{2.8} and \eqref{2.10} and show, by classical arguments, that $\bfv$ is a  solution to \eqref{2.1} that  satisfies, in addition, the ``energy inequality:"
\be
2\|\mathbb D(\bfv)\|_2^2\le\lambda\bftau\cdot\bfe_1\,.
\eeq{2.12}
From \eqref{2.12} we obtain, in particular that $\bftau\neq\0$ because, otherwise, by \lemmref{1.1}, $\bfv\equiv\0$, in contradiction with \eqref{2.1}. Therefore, $\bfv$ is a weak solution to \eqref{1.2} in the sense of \defref{2.1}, and the proof of the theorem is completed.
\par\hfill$\square$
\Br A particular subclass of weak solutions to \eqref{1.2} is the one characterized by having $\bftau=\tau\,\bfe_1$, $\tau>0$, $\bfomega=\0$,  and  $\bfv$ possessing rotational symmetry around $\bfe_1$. The existence of such solutions is shown in \cite[Theorem 1.1]{GaJMP}.
\ER{2.1}\par
The previous theorem proves the existence of a steady-state solution, $\bfv=\bfv(\lambda)$,  to \eqref{1.1} for all $\lambda>0$. However, it is silent about the regularity of the map $\lambda\mapsto \bfv(\lambda)$. In the next theorem, we shall furnish sufficient conditions for the existence  of a local, unique, analytic family of weak solutions  for which the velocity of the center of mass is directed along the same direction. To this end, let $\bfv_\C$ be a weak solution to \eqref{1.2}  corresponding to $\lambda=\lambda_\C$, denote by $\bftau_\C=\tau_\C\bfe_c$ the corresponding translational velocity, where $\bfe_\C:=\bftau_\C/|\bftau_\C|$. Moreover, define
$$
\calv_\C=\calv_\C(\real^3):=\{\bfu\in\calv(\real^3): \ \bfu(x)=\hat{u}_1\bfe_\C+\hat{\bfu}_2\times\bfx\,,\ x\in\Omega_0\}\,,
$$
and set \be\mathscr X_\C(\real^3)=\calv_\C(\real^3)\cap\mathscr X_\sub{\bfe_\C}(\real^3).\eeq{Xc}
We can then prove  that the map $\lambda\mapsto \bfv(\lambda)$ is smooth in a neighborhood of $\lambda_\C$ in the class $\mathscr X_\C$, provided the linearization of \eqref{2.1} around $(\bfv_\C,\lambda_\C)$ is trivial. To this end, let
\be
{\bfscr L}_1:\bfu\in \mathscr X_\C(\real^3)\mapsto \bfscr L_1(\bfu)\in \calv^{-1}(\real^3)\,,
\eeq{L1}
where
\be
\langle \bfscr L_1(\bfu),\bfphi\rangle=2[\bfu, \bfphi]\! -\!\lambda_\C\left[\tau_\C\,\langle \bfe_\C\cdot\nabla\bfu,\bfphi\rangle\!+\!\hat{u}_1\langle\bfe_\C\cdot\nabla\bfv_\C,\bfphi\rangle\! +\!\big(\bfu\cdot\nabla\bfv_\C+\bfv_\C\cdot\nabla\bfu,\bfphi\big)\right]\,,\, \ \bfphi\in\calv(\real^3)\,.
\eeq{L1_1}
The following result holds.
\Bt Let $\bfv_\C$ be a weak solution to \eqref{1.2}  corresponding to $\lambda=\lambda_\C$. Then, the operator $\bfscr L_1$ is Fredholm of index 0. Furthermore, suppose that the equation
\be
\bfscr L_1(\bfu)=\0
\eeq{2.13}
has only the solution $\bfu=\0$ in $\mathscr X_\C(\real^3)$. Then, there exists a neighborhood $U_\C$ of $\lambda_\C$, such that \eqref{1.2} has a  unique family of weak  solutions $\bfv(\lambda)\in\mathscr X_\C(\real^3)$, $\lambda\in U_\C$, which is analytic at $\lambda_\C$ and such that $\bfv(\lambda_\C)=\bfv_\C.$  
\ET{2.2}
{\em Proof.} Consider  the map
$$
\calf :(\bfv,\lambda)\in\mathscr X_\C\times U_\C\mapsto \calf(\bfv,\lambda)\in \calv^{-1}\,,
$$
where
$$
\langle \calf(\bfv,\lambda),\bfphi\rangle:=[\bfv, \bfphi]-\lambda\,\hat{\bfphi}_1\cdot\bfe_1-\lambda\big\langle(\hat{v}_1\bfe_\C-\bfv)\cdot\nabla\bfv,\bfphi\big\rangle\,,\ \ \bfphi\in\calv\,.
$$
The map is well defined. In fact, since $\bfv\in\mathscr X_\C$, we have $\bfe_\C\cdot\nabla\bfv\in\calv^{-1}$. Furthermore, from \eqref{int}, H\"older inequality and \lemmref{1.2} it follows that $\bfv\cdot\nabla\bfv\in \calv^{-1}$ as well. 
In addition, since $\calf$ involves only cubic nonlinearities,  $\calf$ is analytic. We now observe that 
\eqref{2.6} is equivalent to $\calf(\bfv,\lambda)=0$. In \cite[Lemmas 2.2 and 2.3]{GaJMP} it is shown that the operator
$$
\bfscr O:\bfu\in\mathscr X_\C\mapsto \bfscr O(\bfu)\in\calv^{-1}
$$
with 
$$
\langle\bfscr O(\bfu),\bfphi\rangle:= [\bfu, \bfphi] -\lambda_\C\langle\tau_\C\,\bfe_\C\cdot\nabla\bfu,\bfphi\rangle\,,\ \ \bfphi\in\calv\,,
$$
is a homeomorphism, while the operator
$$
\bfscr K:\bfu\in\mathscr X_\C\mapsto \bfscr K(\bfu)\in\calv^{-1}
$$
with 
$$
\langle\bfscr K(\bfu),\bfphi\rangle:=
-\lambda_\C\big\langle(\hat{u}_1\,\bfe_\C-\bfu)\cdot\nabla\bfv_\C-\bfv_\C\cdot\nabla\bfu,\bfphi\big\rangle\,,\ \ \bfphi\in\calv\,,
$$
is compact. As a consequence, $\bfscr L_1\equiv\bfscr O+\bfscr K$ is Fredholm of index 0 and thus, under the assumptions of the theorem, 
it is a homeomorphism. 
By \lemmref{2.1}(d), it is $\calf(\bfv_\C,\lambda_\C)=0$, while the partial Fr\'echet derivative of $\calf$ at $(\bfv_\C,\lambda_\C)$, $D_\sub{\bfv}\calf(\bfv_\C,\lambda_\C)$, is easily shown to satisfy $D_\sub{\bfv}\calf(\bfv_\C,\lambda_\C)=\bfscr L_1$. Since $\bfscr L_1$ is a homeomorphism, the result is an immediate consequence of the analytic version of the Implicit Function Theorem.
\par\hfill$\square$\par
\Br We observe that the assumption of \theoref{2.2} excludes that $(\bfv_\C,\lambda_\C)$ is a steady-state bifurcation point in the class $\mathscr X_\C$ \cite{GaJMP}.
\ER{2.2}
\setcounter{equation}{0}
\section{On the Spectral Properties of the Linearized Operator}
The main objective of this section is to establish some important spectral properties of the  operator obtained by linearizing \eqref{1.2} around the solution $\bfv_\C$. As shown later on, such properties will support one of the basic assumptions of our bifurcation result. To this end, we begin to 
 define the map 
$$
\tilde{\Delta}:\bfw\in W^{2,2}(\Omega)\cap \call^2(\real^3)\mapsto \tilde{\Delta}(\bfw)\in \call^2(\real^3) 
$$
where
\be
\tilde{\Delta}(\bfw)=\left\{\ba{ll}\medskip
 -\Delta\bfw\ \mbox{in $\Omega$}\,,\\ 
\Frac1M\Int{\partial\Omega}{}(2\,\mathbb D(\bfw)\cdot\bfn)+\left(\Frac1\cali\Int{\partial\Omega}{}\bfy\times(2\,\mathbb D(\bfw)\cdot\bfn)\right)\times\bfx\ \ \mbox{in $\Omega_0$\,,}\ea\right.
\eeq{Delta}
and  set $\bfscr A:=\mathscr P\,\tilde{\Delta}$, with $\mathscr P$  the orthogonal projection of $\call^2(\real^3)$ onto $\calh(\real^3)$; see \lemmref{0.1}. We next consider the operator 
\be 
\bfscr L_0:\bfu\in Z^{2,2}:=W^{2.2}(\Omega)\cap\calh(\real^3)\mapsto\call_0(\bfu)= -\lambda_\C\bftau_\C\cdot\nabla\bfu+\bfscr A(\bfu)\in \calh(\real^3)\,.
\eeq{L0}
$\mathscr L_0$ is well defined since \be\bftau_\C\cdot\nabla\bfu\in \calh(\real^3).\eeq{tau}
In fact, we observe that $\bfu\in Z^{2,2}$ implies $\bfu\in \calv(\real^3)$, so that 
by \lemmref{0.1} and \eqref{0.0}, it follows that \eqref{tau} reduces to prove that
\be
\int_\Omega\bftau_\C\cdot\nabla\bfu\cdot\bfh+M\,\widehat{(\bftau_\C\cdot\nabla\bfu)}_1\cdot\hat{\bfh}_1+\cali\,\widehat{(\bftau_\C\cdot\nabla\bfu)}_2\cdot\hat{\bfh}_2=0\,,\ \mbox{for all $\bfh\in\calg(\real^3)$}\,.
\eeq{dec0}
Taking into account that in $\Omega_0$ it is $\bftau_\C\cdot\nabla\bfu=\bftau_\C\cdot\nabla (\hat{\bfu}_2\times\bfx)=\hat{\bfu}_2\times\bftau_\C$, while, by \eqref{spazi}$_3$, $\bfh=\nabla p$ in $\Omega$ with $\hat{\bfh}_1=-M^{-1}\int_{\partial\Omega}p\,\bfn$ and $\hat{\bfh}_2= -\cali^{-1}\int_{\partial\Omega}p\,\bfx\times\bfn$,  \eqref{dec0} becomes
$$
\int_\Omega\bftau_\C\cdot\nabla\bfu\cdot\nabla p -\int_{\partial\Omega}p\,\hat{\bfu}_2\times\bftau_\C\cdot\bfn=0\,,
$$
whose validity is immediately checked by integrating by parts the volume integral, and recalling that $\Div\bfu=0$.
\par
The following preliminary result holds.
\Bl Let $\zeta\in\real\backslash\{0\}$ and $(\lambda_\C,\bftau_\C)\in\real\times\real^3$. Then, the operator $\mathscr L_0+\i\,\zeta$ is a homeomorphism of $Z_{\mathbb C}^{2,2}$ onto $\calh_{\mathbb C}(\real^3)$.  Moreover, there is $c=c(\Omega,\rho_{\mathscr S}/\rho_{\mathscr L})$, such that
\be
\|D^2\bfu\|_2+|\zeta|^{\frac12}\|\nabla\bfu\|_2+|\zeta|(\|\bfu\|_2+M^{\frac12}|\bfchi|+\cali^{\frac12}|\bfsigma|)\le c\,\|(\bfscr L_0+\i\,\zeta)(\bfu)\|_2\,,\ \ |\zeta|\ge \max\{\lambda_\C^2 |\bftau_\C|^2,1\}\,,
\eeq{4.2}
where $\bfchi:=\hat{\bfu}_1$, $\bfsigma:=\hat{\bfu}_2$.
\EL{4.1}
{\em Proof.} The homeomorphism property can be obtained by proving  that for any $\bfcalf\in\calh_{\mathbb C}(\real^3)$, there exists one and only one $\bfu\in Z_{\mathbb C}^{2,2}$ such that $\mathscr L_0(\bfu)+\i\,\zeta\,\bfu=\bfcalf$. In view of \lemmref{aria}, \lemmref{0.1} and \eqref{Delta} this is equivalent to requiring that
for any $(\bff,\bfF,\bfG)\in L^2_{\mathbb C}(\Omega)\times \mathbb C\times\mathbb C$, with $\Div\bff=0$,
the problem 
\be\ba{cc}\medskip\left.\ba{ll}\medskip\Delta\bfu+\lambda_\C\bftau_\C\cdot\nabla\bfu-\nabla{\mathfrak p}=\i\,\zeta\,\bfu +\bff\\
\Div\bfu=0\ea\right\}\,  \mbox{in $\Omega$}\,,\\ \medskip
\bfu=\bfchi+\bfsigma\times\bfx\ \ \mbox{at $\partial\Omega$}\,,\\ 
\i\,\zeta\,M\,{\bfchi}+\Int{\partial\Omega}{}\mathbb T(\bfu,{\mathfrak p})\cdot\bfn=\bfF\,,\ \   
\i\,\zeta\,\mathcal I\,{\bfsigma}+\Int{\partial\Omega}{}\bfx\times\mathbb T(\bfu,{\mathfrak p})\cdot\bfn=\bfG\,,
\ea
\eeq{4.1}
has one and only one solution $(\bfu,\mathfrak p,\bfchi,\bfsigma)\in W_{\mathbb C}^{2,2}(\Omega)\times D^{1,2}_{\mathbb C}(\Omega)\times\mathbb C\times\mathbb C$.
Let us dot-multiply both sides of \eqref{4.1}$_1$ by $\bfu^*$ (${}^*=$\,complex conjugate) and integrate by parts over $\Omega$. Taking into account \eqref{4.1}$_{2-5}$ we show
\be
2\|\mathbb D(\bfu)\|_2^2+\i\,\zeta\,(\|\bfu\|_2^2+M|\bfchi|^2+\cali|\bfsigma|^2)-\lambda(\bftau\cdot\nabla\bfu,\bfu^*)=(\bff,\bfu^*)+\bfF\cdot\bfchi^*+\bfG\cdot\bfsigma^*\,,
\eeq{4.3}
where, here and in what follows, in order to simplify the notation we suppress the subscript $\C$. Taking the real and imaginary parts of \eqref{4.2},  using \eqref{1.8} and Schwarz inequality, and observing that $\Re \,(\bftau\cdot\nabla\bfu,\bfu^*)=0$, we deduce
\be
\|\nabla\bfu\|_2^2\le \bfpzc F\,\bfpzc W\,;\ \ \
|\zeta|\,\|\bfu\|^2\le (\lambda\,|\bftau|\,\|\nabla\bfu\|_2+\bfpzc F)\,\bfpzc W\,,
\eeq{4.4}
where
$$\ba{ll}\medskip
\|\bfW\|:=\big(\|\bfu\|_2^2+M|\bfchi|^2+\cali|\bfsigma|^2\big)^{\frac12}\\ \medskip
\bfpzc W:=\|\bfu\|_2+M^{\frac12}|\bfchi|+\cali^{\frac12}|\bfsigma|\,,\\
\bfpzc F:=\|\bff\|_2+M^{-\frac12}|\bfF|+\cali^{-\frac12}|\bfG|\,.
\ea
$$
Observing that 
$\bfpzc W\le 3^{\frac12}\|\bfW\|$, from \eqref{4.4} and Cauchy-Schwarz inequality we get
$$
|\zeta|\,\|\bfu\|\le 3^{\frac12}(\lambda\,|\bftau|\,\bfpzc F^{\frac12}\,\bfpzc W^{\frac12}+\bfpzc F)\le 3^{\frac12}\left(\frac{3\lambda^2\,|\bftau|^2}{|\zeta|}+1\right)\bfpzc F+\half|\zeta|\,\|\bfW\|\,,
$$
from which we conclude
\be
|\zeta|\,\|\bfW\|\le  3^{\frac12}\left(\frac{3\lambda^2\,|\bftau|^2}{|\zeta|}+2\right)\bfpzc F\,.
\eeq{4.5}
Replacing the latter into \eqref{4.4}$_1$, we infer
\be
|\zeta|^{\frac12}\|\nabla\bfu\|_2\le 3^{\frac14} \left(\frac{3\lambda^2\,|\bftau|^2}{|\zeta|}+2\right)^{\frac12}\bfpzc F\,.
\eeq{4.6}
Combining the estimates \eqref{4.5}, \eqref{4.6} with classical Galerkin method, we can proceed as in the proof of \theoref{2.1} and show that for any given $(\bff,\bfF,\bfG)$ in the specified class and $\zeta\neq 0$, there exists a (unique, weak) solution to \eqref{4.1} such that $\bfu\in \calv_{\mathbb C}(\Omega)\cap L^2_{\mathbb C}(\Omega)$, satisfying  \eqref{4.5}, \eqref{4.6}. We next write \eqref{4.1}$_{1-3}$ as the following Stokes problem
$$\ba{cc}\medskip\left.\ba{ll}\medskip\Delta\bfu=\nabla{\mathfrak p}+\bfcalg\\
\Div\bfu=0\ea\right\}\,  \mbox{in $\Omega$}\,,\\ \medskip
\bfu=\bfchi+\bfsigma\times\bfx\ \ \mbox{at $\partial\Omega$}
\ea
$$
where
$$
\bfcalg:=-\lambda\, \bftau\cdot\nabla\bfu +\i\,\zeta\,\bfu+\bff\,.
$$
Since $\bfcalg\in L^2_{\mathbb C}(\Omega)$ and $\bfu\in W^{1,2}_{\mathbb C}(\Omega)$, from classical results \cite[Theorems IV.5.1 and V.5.3]{Gab} it follows that $D^2\bfu\in L^2(\Omega)$, thus completing the existence (and uniqueness) proof. Furthermore, by \cite[Lemma IV.1.1 and V.4.3]{Gab} we get  
\be
\|D^2\bfu\|_2\le c\,\big[\|\bff\|_2+(\lambda\,|\bftau|+1)\|\nabla\bfu\|_2+(|\zeta|+1)\,\|\bfu\|_2+|\bfchi|+|\bfsigma|\big]\,.
\eeq{4.7}
As a result, if $|\zeta|\ge \max\{
\lambda^2|\bftau|^2,1\}$, the inequality in \eqref{4.2} is a consequence of \eqref{4.5}--\eqref{4.7}.
\par\hfill$\square$\par
Let 
\be
\bfcalk: \bfu\in  Z^{2,2} \mapsto \bfcalk(\bfu)\in \call^2(\real^3)
\eeq{3.32}
where ($\bfchi:=\hat{\bfu}_1$), 
\be
\bfcalk(\bfu)=\left\{\ba{ll}\medskip
\lambda_\C\big(\bfv_\C\cdot\nabla\bfu+(\bfu-\bfchi)\cdot\nabla\bfv_\C\big)\ \ \mbox{in $\Omega$}\,,\\
\0\ \ \mbox{in $\Omega_0$}\,.
\ea\right.
\eeq{3.32_1}
From  \lemmref{2.1}, we get, in particular,
\be
\bfv_\C\in L^\infty(\Omega)\cap L^4(\Omega)\,,\ \ \nabla\bfv_\C\in L^\infty(\Omega)\cap D^{1,2}(\Omega)\,, 
\eeq{stoc}
and so we easily check that the operator $\bfcalk$ is well defined.
Finally, let
\be
\bfscr L_2:\bfu\in W^{2,2}(\Omega)\cap\calh(\real^3)\mapsto \bfscr L_2(\bfu)=-\lambda_\C\bftau_\C\cdot\nabla\bfu+\bfscr A(\bfu)+\mathscr P\,\bfcalk(\bfu)\in \calh(\real^3)
\eeq{L2}
We are now ready to show the main result of this section.
\Bt The operator
$ 
\bfscr L_2+\i\,\zeta
$ 
is Fredholm of index 0, for all $\zeta\neq0$. Moreover, let $\sigma(\bfscr L_2)$ be the spectrum of $\bfscr L_2$. Then,  $\sigma(\bfscr L_2)\cap \{\i\,\real\backslash\{0\}\}$
consists, at most, of a finite or countable number of eigenvalues, each of which is isolated and of finite (algebraic) multiplicity, that can
only accumulate at 0.  
\ET{5.1}
{\em Proof.} We show that the operator $\bfcalk$ defined in \eqref{3.32}, \eqref{3.32_1} is compact. Let $\{\bfu_k\}$ be a sequence bounded in $Z^{2,2}$. This implies, in particular, that there is $M>0$  independent of $k$ such that ($\bfchi_k:=\widehat{(\bfu_k)}_1$)
\be
\|\bfu_k\|_{2,2}+|\bfchi_k|\le M\,.
\eeq{5.9}
The latter, along with the compact embedding $W^{2,2}(\Omega)\subset W^{1,2}(\Omega_R)$, for all $R>1$ entails the existence of $(\bfu_*,\bfchi_*)\in W^{2,2}(\Omega)\times\real$ and subsequences, again denoted by $\{\bfu_k,\bfchi_k)$ such that
\be
\bfu_k\to\bfu_*\ \ \mbox{strongly in $W^{1,2}(\Omega_R)$, for all $R>1$}\,;\ \ \bfchi_k\to\bfchi_*\ \ \mbox{in $\real^3$}\,. 
\eeq{JN}
Without loss, we assume $\bfu_*\equiv\bfchi_*\equiv\0$. By H\"older inequality and \eqref{3.32}, we deduce 
$$
\|\bfcalk(\bfu_k)\|_2\le \lambda_\C\left[(\|\bfv_\C\|_\infty+\|\nabla\bfv_\C\|_\infty)(\|\bfu_k\|_{1,2,\Omega_R}+|\bfchi_k|)+\|\bfv_\C\|_{1,4,\Omega^R}\|\nabla\bfu_k\|_{1,4}\right]
$$
Therefore, passing to the limit $k\to\infty$ in the previous inequality, and using \eqref{stoc}, \eqref{JN}, \eqref{5.9} and the embedding $W^{2,2}\subset W^{1,4}$,  we infer
$$
\lim_{k\to\infty}\|\bfcalk(\bfu_k)\|_2\le M\,\|\bfv_\C\|_{1,4,\Omega^R}\,,
$$
which, in turn, again by \eqref{stoc} and the arbitrariness of $R>1$ shows the desired property. From this and \lemmref{4.1} it then follows that the operator
$$ 
\widehat{\bfscr L}_\zeta:=\bfscr L_1+\i\,\zeta
$$ 
is Fredholm of index 0, for all $\zeta\neq0$. The theorem is then a consequence of well-known results (e.g. \cite[Theorem XVII.4.3]{GG}) provided we show that the null space of $\widehat{\bfscr L}_\zeta$ is trivial, for all sufficiently large $|\zeta|$. To this end, we observe that $\widehat{\bfscr L}_\zeta(\bfu)=0$ means $\bfscr L_0(\bfu)+\i\,\zeta\,\bfu=-\bfcalk(\bfu)$. Applying \eqref{4.2}, we thus get, in particular, the following inequality valid for all sufficiently large $|\zeta|$
$$
|\zeta|^{\frac12}\|\nabla\bfu\|_2+|\zeta|(\,\|\bfu\|_2+|\bfchi|)\le c\,\|\bfcalk(\bfu)\|_2\,,
$$ 
where $c$ is independent of $\zeta$. Also, from \eqref{3.32}, \eqref{3.32_1}, \eqref{stoc} and H\"older inequality, we have
$$
\|\bfcalk(\bfu)\|_2\le \lambda_\C\, \|\bfv_\C\|_{1,\infty}(\|\bfu\|_{1,2}+|\bfchi|)\,,
$$
and so, from the last two displayed equations we deduce $\bfu\equiv\0$, provided we choose $|\zeta|$ larger than a suitable positive constant depending only on $\Omega, \lambda_\C, M$, and $\cali$. This completes the proof of the theorem.
\par\hfill$\square$\par

\setcounter{equation}{0}
\section{On the Linearized Time-Periodic Operator}
Let
$$\ba{ll}\medskip
\mathfrak L^2_\sharp:=\{\bfcalf\in L^2(\call^2)\,, \ \bfcalf\ \mbox{is $2\pi$-periodic with}\ \bar{\bfcalf|_{\Omega}}=\bar{\widehat{\bfcalf}_1}=\bar{\widehat{\bfcalf}_2}=\0\}\,,\\ \medskip
\mathfrak H_\sharp:=\{\bfpzc f\in\mathfrak L^2_\sharp:\ \bfpzc f\in L^2(\calh)\}\\
\mathfrak W_\sharp^2:=\{\bfw\in \mathcal W_\sharp^2: \bfu\in W^{1,2}(\calh)\}\,,
\ea
$$
and, for $\zeta_0\neq 0$, consider the operators
$$
\bfcalq_0:\bfw\in \mathfrak W_\sharp^2\mapsto \zeta_0\,\partial_\s\bfw+\bfscr L_0(\bfw)\in \mathfrak H_\sharp\,,
$$
and
\be
\bfcalq:\bfw\in \mathfrak W_\sharp^2\mapsto \zeta_0\,\partial_\s\bfw+\bfscr L_2(\bfw)\in \mathfrak H_\sharp\,,
\eeq{Q}
where $\mathscr L_0$ and $\mathscr L_2$ are given in \eqref{L0} and \eqref{L2}, respectively.
\par
Our objective in this section is to establish some important properties for both operators. 
We begin to show the following lemma.
\Bl Let $\tau\in\real$. Then, the boundary-value problems, $i\in\{1,2,3\}$, $k\in \mathbb Z\backslash\{0\}$,
\be
\ba{cc}\medskip\left.\ba{ll}\medskip
{\rm i}\,k\,\bfh_k^{(i)}-\tau\partial_1\bfh_k^{(i)}=\Delta \bfh_k^{(i)}-\nabla p_k^{(i)}\\
\Div \bfh_k^{(i)}=0\ea\right\}\ \ \mbox{in $\Omega$}\\
\bfh_k^{(i)}=\bfe_i\ \ \mbox{at $\partial\Omega$}\,,\ \ \bfh_0^{(i)}=\0
\ea
\eeq{3.2}
and
\be
\ba{cc}\medskip\left.\ba{ll}\medskip
{\rm i}\,k\,\bfH_k^{(i)}-\tau\partial_1\bfH_k^{(i)}=\Delta \bfH_k^{(i)}-\nabla P_k^{(i)}\\
\Div \bfH_k^{(i)}=0\ea\right\}\ \ \mbox{in $\Omega$}\\
\bfH_k^{(i)}=\bfe_i\times\bfx\ \ \mbox{at $\partial\Omega$}\,,\ \ \bfH_0^{(i)}=\0
\ea
\eeq{3.3}
have unique solutions $(\bfh_k^{(i)},p_k^{(i)}),(\bfH_k^{(i)},P_k^{(i)})\in W^{2,2}(\Omega)\times D^{1,2}(\Omega)$. These solutions satisfy the estimates
\be\ba{ll}\medskip
\|\bfh_k^{(i)}\|_2+\|\bfH_k^{(i)}\|_2\le C\\ \medskip
\|\nabla\bfh_k^{(i)}\|_2+\|\nabla\bfH_k^{(i)}\|_2\le C \,(|k|+1)^{\frac12}\\
\|D^2\bfh_k^{(i)}\|_2+\|D^2\bfH_k^{(i)}\|_2\le C \,(|k|+1)\,,
\ea
\eeq{3.4}
where $C$ is a constant independent of $k$. Moreover, for fixed $k$, consider the $3\times3$ matrices $\mathbb K,$ $\mathbb A$, $\mathbb P$ and $\mathbb S$ defined by the components ($j,i=1,2,3$):
\be\ba{cc}\ms
({\mathbb K})_{ji}=\IdS(\mathbb T(\bfh^{(i)}_k,p^{(i)}_k)\cdot\bfn)_j,\\ \ms
({\mathbb A})_{ji}=\IdS(\bfx\times\mathbb T(\bfH^{(i)}_k,P^{(i)}_k)\cdot\bfn)_j\\ \ms
({\mathbb P})_{ji}=\IdS(\bfx\times\mathbb T(\bfh^{(i)}_k,p^{(i)}_k)\cdot\bfn)_j\\ ({\mathbb S})_{ji}=\IdS(\mathbb T(\bfH^{(i)}_k,P^{(i)}_k)\cdot\bfn)_j
\ea
\eeq{Matrices}
and define the $6\times 6$ matrix $\bfpzc A$ as follows
$$
\bfpzc A:=\left(\ba{cc}\medskip \mathbb K\ \ \mathbb P\\
\mathbb S\ \ \mathbb A\ea\right)\,.
$$
Then, for any $\mu\in\real$, both $\mathbb K+\i\,\mu\mathds{1}$ and $\mathbb A+\i\,\mu\mathds{1}$ are invertible. Moreover,  for every $\bfzeta\in \mathbb C^6$, we have
\be
{\rm i}\,k\,\|\bfsf v\|_2^2+2\|\mathbb D(\bfsf v)\|_2^2-\tau(\partial_1\bfsf v,\bfsf v^*)=\bfzeta^*\cdot\bfpzc A\cdot\bfzeta
\eeq{3.A}
where $\bfsf v:= \zeta_i\bfh^{(i)}_k+\zeta_{i+3}\bfH^{(i)}_k$. Finally, for every $(\lambda,\mu)\in\real\times\real$, the matrix
$$
\bfpzc A+\left(\ba{cc}\medskip \i\lambda\,\mathbb I & 0\\
0& \i\mu\mathbb I\ea\right):=\bfpzc A+\bfpzc J
$$
is invertible.
\EL{3.1}
{\em Proof.} We begin to show the estimate for $\bfh^{(i)}_k$. Since the proof is the same
for $i = 1, 2, 3$, we chose $i = 1$ and, for simplicity, omit the superscript. Let
$\phi = \phi(|x|)$ be a (smooth) cut--off function such that
$$
\phi(|x|)=\left\{\ba{ll}\medskip 1\ \mbox{in $\Omega_R$}\\
0\ \mbox{in $\bar{\Omega^{2R}}$}\ea\right.
$$
and set $\bfPhi(x)=\curl\big(x_2\phi(|x|)\bfe_3\big)$. Clearly, $\Div\bfPhi=0$ and $\bfPhi(x)=\bfe_1$ in a neighborhood of $\partial\Omega$. Moreover $\bfPhi(x)\equiv\0$ in $\Omega^{2R}$. Setting $\bfv_k:=\bfh_k-\bfPhi$, from \eqref{3.2} we deduce that $\bfv_k$ solves the following boundary-value problem, for all $|k|\ge1$:
\be
\ba{cc}\medskip\left.\ba{ll}\medskip
{\rm i}\,k\,\bfv_k-\tau\,\partial_1\bfv_k=\Delta \bfv_k-\nabla p_k^{(i)}+\tau\partial_1\bfPhi-{\rm i}\,k\,\bfPhi+\Delta\bfPhi\\
\Div \bfv_k=0\ea\right\}\ \ \mbox{in $\Omega$}\\
\bfv_k=\0\ \ \mbox{at $\partial\Omega$}\,.
\ea
\eeq{3.5}
Existence to \eqref{3.5} in the stated function class can be easily obtained by the
Galerkin method combined with the estimate that we are about to derive.
Let us dot-multiply both sides of \eqref{3.5}$_1$ by $\bfv^*$ where the star denotes c.c. After integrating by parts as necessary, we get
\be
{\rm i}\,k\,\|\bfv_k\|_2^2-\tau(\partial_1\bfv,\bfv_k^*)+\|\nabla\bfv_k\|_2^2=(\bfcalf_k,\bfv_k^*)\,,
\eeq{3.6}
where $\bfcalf_k:=\tau\partial_1\bfPhi-{\rm i}\,k\,\bfPhi+\Delta\bfPhi$. We next observe that, by the properties of $\bfPhi$,
\be
\|\bfcalf_k\|_2\le c\,(|k|+1)
\eeq{3.7}
where, here and in the rest of the proof, $c$ denotes a generic (positive) constant independent of $k$. Also, by means of an integration by parts, we show
\be
\Re (\partial_1\bfv_k,\bfv_k^*)=0\,.
\eeq{3.8}
Thus, by taking the real part of \eqref{3.6} and using \eqref{3.7} and \eqref{3.8} we infer
\be
\|\nabla\bfv_k\|_2^2\le c\,(|k|+1)\|\bfv_k\|_2\,.
\eeq{3.9}
Likewise, taking the imaginary part of \eqref{3.6} and employing \eqref{3.7}--\eqref{3.9} along with Schwarz inequality, we obtain 
$$
|k|\|\bfv\|_2\le c\,(\|\nabla\bfv_k\|_2+|k|+1)\le c(|k|+1)\|\bfv_k\|_2^{\frac12}\,,
$$
which implies
\be
\|\bfv\|_2\le c\,.
\eeq{3.10}
Taking into account that $\bfh_k=\bfPhi+\bfv_k$, \eqref{3.10} proves \eqref{3.4}$_1$ for $\bfh_k$. Similarly, replacing \eqref{3.10} into \eqref{3.9}, we arrive at \eqref{3.4}$_2$. Finally, from classical estimates on the Stokes problem \cite[Lemma 1]{Hey} we find
$$
\|D^2\bfv_k\|_2\le c\,\big(\|\partial_1\bfv\|_2+\|\bfcalf_k\|_2+\|\nabla\bfv\|_2\big)
$$
and so \eqref{3.4}$_3$ follows from this inequality, \eqref{3.7}, \eqref{3.9} and \eqref{3.10}. Concerning the fields $\bfH^{(i)}_k$, let $\bfPsi^{(i)}=\phi(|x|)\bfe_i\times\bfx$, and set $\bfV^{(i)}_k:=\bfH^{(i)}+\bfPsi^{(i)}$. Obviously, the support of $\bfPsi^{(i)}$ is contained in $\Omega_{2R}$, $\Div\bfPsi^{(i)}=0$ and $\bfPsi^{(i)}|_{\partial\Omega}=\bfe_i\times\bfx$. Thus, from \eqref{3.3} it follows that 
$\bfV^{(i)}_k$ is a solution to \eqref{3.5} with $(\bfV^{(i)}_k,\bfPsi^{(i)})$ in place of $(\bfv_k,\bfPhi)$. Therefore, we can use exactly the same arguments used earlier in the proof to show that also $\bfH^{(i)}_k$ satisfies the stated properties. Let $\bfalpha\in\mathbb C^3$, and, for fixed $k\neq0$,
 set\footnote{Summation over repeated indices.}
$$
\bfsf u:= \alpha_i\bfh^{(i)}_k\,,\ \ {\sf q}:=\alpha_i\,p^{(i)}_k\,.
$$
From \eqref{3.2}  we then find
\be
\ba{cc}\medskip\left.\ba{ll}\medskip
{\rm i}\,k\,\bfsf u-\tau\partial_1\bfsf u=\Div\mathbb T(\bfsf u,{\sf q})\\
\Div \bfsf u=0\ea\right\}\ \ \mbox{in $\Omega$}\\
\bfsf u=\bfalpha\,\ \mbox{at $\partial\Omega$}.
\ea
\eeq{3.2_ar}
Dot-multiplying both sides of \eqref{3.2_ar}$_1$ by $\bfsf u^*$ and integrating by parts over $\Omega$ we deduce
$$
\i\,k\|\bfsf u\|_2^2+\|\mathbb D(\bfsf u)\|_2^2-\tau(\partial_1\bfsf u,\bfsf u^*)=\bfalpha^*\cdot\mathbb K\cdot\bfalpha\,.
$$
Now, suppose that there is $\hat{\bfalpha}\in\mathbb C$ such that $\mathbb K\cdot\hat{\bfalpha}=-\i\,\mu\,\hat{\bfalpha}$, for some $\mu\in\real$. Then from the previous relation we obtain
$$
\i\left(k\|\bfsf u\|_2^2+\mu\,|\hat{\bfalpha}|^2\right)-\tau(\partial_1\bfsf u,\bfsf u^*)=\|\mathbb D(\bfsf u)\|_2^2\,,
$$
which, in turn, taking into account that $\Re\,(\partial_1\bfsf u,\bfsf u^*)=0$, allows us to we deduce $\bfsf u=\0$ in $W^{2,2}(\Omega)$. The latter implies $\hat{\bfalpha}=\0$ and thus shows the desired property for $\mathbb K$.  In a similar manner, we prove the same property for $\mathbb S$. 
Next, let $\bfzeta\in \mathbb C^6$ and define\footnote{Summation over repeated indices.}
\be
\bfsf v:= \zeta_i\bfh^{(i)}_k+\zeta_{i+3}\bfH^{(i)}_k\,,\ \ {\sf p}:=\zeta_i\,p^{(i)}_k+\zeta_{i+3}\,P^{(i)}_k\,.
\eeq{3.12}
Employing \eqref{3.2} and \eqref{3.3} we then deduce
\be
\ba{cc}\medskip\left.\ba{ll}\medskip
{\rm i}\,k\,\bfsf v-\tau\partial_1\bfsf v=\Div\mathbb T(\bfsf v,{\sf p})\\
\Div \bfsf v=0\ea\right\}\ \ \mbox{in $\Omega$}\\
\bfsf v=\zeta_i\bfe_i+(\zeta_{3+i}\bfe_i)\times\bfx\ \ \mbox{at $\partial\Omega$}\,.
\ea
\eeq{3.13}
By dot-multiplying both sides of \eqref{3.13}$_1$ by $\bfsf v^*$ and integrating by parts over $\Omega$, we show
$$
{\rm i}\,k\,\|\bfsf v\|_2^2+2\|\mathbb D(\bfsf v)\|_2^2-\tau(\partial_1\bfsf v,\bfsf v^*)=\int_{\partial\Omega}\bfsf v^*\cdot\mathbb T(\bfsf v,{\sf p})\cdot\bfn\,.
$$
As a consequence, \eqref{3.A} follows by replacing \eqref{3.12} in this last inequality and using \eqref{Matrices} and \eqref{3.13}$_3$. Finally, let $\bfzeta\in\mathbb C^6$ such that $(\bfpzc A+\bfpzc J)\cdot\bfzeta=\0$. From \eqref{3.A} it then follows that 
$$
{\rm i}\big[\,k\,\|\bfsf v\|_2^2+\sum_{i=1}^3(\lambda|\zeta_i|^2+\mu|\zeta_{i+3}|^2)\big] +2\|\mathbb D(\bfsf v)\|_2^2-\tau(\partial_1\bfsf v,\bfsf v^*)=0\,,
$$ 
from which, recalling that $\Re(\partial_1\bfsf v,\bfsf v^*)=0$, and that  $\bfsf v$ is a solution to \eqref{3.13} in the space $W^{1,2}(\Omega)$, we at once obtain $\bfsf v\equiv \0$, implying $\bfzeta=\0$.  The proof is completed.
\par\hfill$\square$\par
\Br Even though the results of the previous lemma are stated for $\Omega$ the exterior of a ball, the reader will check with no effort that they continue to hold --without changes in their proof-- for any exterior domain of class $C^2$. Therefore, they generalize those obtained in \cite[Lemma 5.1]{GaSP}
\ER{3.1}\par
With the help of \lemmref{3.1}, we are now able to show the following one.
\Bl Let $\tau\in\real$. Then, for any $(\bff,\bfF,\bfG)\in \call_\sharp^{2}\times L^2_\sharp\times L^2_\sharp$, the  problem 
\be\ba{cc}\medskip\left.\ba{ll}\medskip\partial_\s\bfw-\bftau\cdot\nabla\bfw=\Delta\bfw-\nabla{\phi}+\bff\\
\Div\bfw=0\ea\right\}\,  \mbox{in $\Omega\times [0,2\pi]$}\,,\\ \medskip
\bfw=\bfchi+\bfsigma\times\bfx\ \ \mbox{at $\partial\Omega\times [0,2\pi]$}\,,\\ 
\,M\dot{\bfchi}+\Int{\partial\Omega}{}\mathbb T(\bfw,{\phi})\cdot\bfn=\bfF\,,\ \   
\,\mathcal I\,\dot{\bfsigma}+\Int{\partial\Omega}{}\bfx\times\mathbb T(\bfw,{\phi})\cdot\bfn=\bfG\, \ \ \mbox{in $[0,2\pi]$}\,,
\ea
\eeq{3.1}
has one and only one solution $\big(\bfw,\phi,\bfchi,\bfsigma\big)\in \mathcal W_\sharp^{2}
\times \mathcal P^{1,2}\times W^{1,2}_\sharp\times W^{1,2}_\sharp$. This solution satisfies the estimate
\be
\|\bfw\|_{\mathcal W_\sharp^{2}}+\|{\phi}\|_{\mathcal P^{1,2}}+\|\bfchi\|_{W^{1,2}}+\|\bfsigma\|_{W^{1,2}}\le C_2\,\Big(\|\bff\|_{\call_\sharp^{2}}+\|\bfF\|_{L^2}+\|\bfG\|_{L^2}\Big)\,,
\eeq{n}
where $C_2=C_2(\Omega,\tau,\rho_{\mathscr S}/\rho_{\mathscr L})$.
\EL{3.2}
{\em Proof.} Since the actual values of  $M$ and $\cali$ are irrelevant to the proof, we put, for simplicity, $M=\cali=1$.  Moreover, without loss of generality, we may take $\bftau=\tau\,\bfe_1$. Let $\bfw=\bfz+\bfu$ where $\bfz$ and $\bfu$ satisfy the following set of equations
\be\ba{cc}\medskip\left.\ba{ll}\medskip
\partial_\s\bfz-\tau\partial_1\bfz-\Delta\bfz=-\nabla {\sf r}+\bff\\ 
\Div\bfz=0\ea\right\}\ \ \mbox{in $\Omega\times [0,2\pi]$}\\
\bfz|_{\partial\Omega}=\0
\ea
\eeq{2.3} 
and
\be\ba{cc}\medskip\left.\ba{ll}\medskip
\partial_\s\bfu-\tau\partial_1\bfu-\Delta\bfu=-\nabla {\sf q}\\ 
\Div\bfu=0\ea\right\}\ \ \mbox{in $\Omega\times[0,2\pi]$}\\ \medskip
\bfu|_{\partial\Omega}=\bfchi+\bfsigma\times\bfx\,;\\ \medskip\dot{\bfchi}+\Int{\partial\Omega}{}\bfT(\bfu,{\sf q})\cdot\bfn=\bfF-\Int{\partial\Omega}{}\bfT(\bfz,{\sf r})\cdot\bfn:=\bfcalf\,,\ \ \mbox{in $[0,2\pi]$}\,,\\
\dot{\bfsigma}+\Int{\partial\Omega}{}\bfx\times\mathbb T(\bfu,{\sf q})\cdot\bfn=\bfG-\Int{\partial\Omega}{}\bfx\times\mathbb T(\bfz,{\sf r})\cdot\bfn:=\bfcalg\,,\ \ \mbox{in $[0,2\pi]$}\,.
\ea
\eeq{3.17}
From \cite[Theorem 12]{GaMaH}, it follows that there exists a unique solution $(\bfz,\tau)\in  \mathcal W_\sharp^{2}\times\calp^{1}$   that, in addition, obeys the inequality
\be
\|\bfz\|_{\mathcal W_\sharp^{2}}+\|\sf r\|_{\mathcal P^{1,2}}\le c\,\|\bff\|_{\mathcal L_\sharp^{2}}\,.
\eeq{3.18}
Furthermore, by trace theorem\Footnote{Possibly, by modifying $\sf r$ by adding to it a suitable function of time.} and \eqref{3.18} we get
$$
\|\Int{\partial\Omega}{}\mathbb T(\bfz,{\sf \sf r})\cdot\bfn\|_{L^2}+\|\Int{\partial\Omega}{}\bfx\times \mathbb T(\bfz,{\sf \tau})\cdot\bfn\|_{L^2}\le c\,\left(\|\bfz\|_{\mathcal W^{2}_\sharp}+\|{\sf \sf r}\|_{\mathcal P^{1,2}}\right)\le c\,\|\bff\|_{\mathcal L_\sharp^{2}}\,,
$$
so that both functions $\bfcalf$ and $\bfcalg$ in \eqref{3.17} are in $L^2_\sharp(0,2\pi)$ and satisfy
\be
\|\bfcalf\|_{L^2}+\|\bfcalg\|_{L^2}\le c(\|\bff\|_{\call^2_\sharp}+\|\bfF\|_{L^2}+\|\bfG\|_{L^2})
\,.
\eeq{3.19} 
To find solutions to \eqref{3.17}, we
formally expand $\bfu,$ ${\sf q}$,   $\bfchi$ and $\bfsigma$ in Fourier series (summation over repeated indices):
\be\ba{ll}\medskip
\bfu(x,\s)=\bfu_k(x)\,{\rm e}^{\i k\,\s}\,,\ \ {\sf q}(x,\s)={\sf q}_k(x)\,{\rm e}^{\i k\,\s}\,,\ \
\bfchi(\s)=\bfchi_k \,{\rm e}^{\i k\,\s}\,,\ \ \bfsigma(\s)=\bfsigma_k \,{\rm e}^{\i k\,\s}\,,\ k\in\mathbb Z\backslash\{0\}\,,\\ \bfu_0\equiv\nabla{\sf q}_0\equiv\bfchi_0\equiv\bfsigma_0\equiv\0\,,
\ea
\eeq{3.20}
where $(\bfu_k,{\sf q}_k,\bfchi_k,\bfsigma_k)$ solve the problem
($k\neq0$)
\be\ba{cc}\medskip\left.\ba{ll}\medskip
\i\,k\,\bfu_k-\tau\partial_1\bfu_k=\Delta \bfu_k-\nabla{\sf q}_k\\
\Div\bfu_k=0\ea\right\}\ \ \mbox{in $\Omega$}\\
\bfu_k|_{\partial\Omega}=\bfchi_k+\bfsigma_k\times\bfx\,,
\ea
\eeq{3.21}
subject to the further conditions
\be
\i\,k\,{\bfchi_k}+\Int{\partial\Omega}{}\mathbb T(\bfu_k,{\sf q}_k)\cdot\bfn=\bfcalf_k\,,\ \ \i\,k\,{\bfsigma_k}+\Int{\partial\Omega}{}\bfx\times\mathbb T(\bfu_k,{\sf q}_k)\cdot\bfn=\bfcalg_k\,,
\eeq{3.22}
with $\{\bfcalf_k\},\{\bfcalg_k\}$ are Fourier coefficients of $\bfcalf$ and $\bfcalg$, respectively, and $\bfcalf_0\equiv\bfcalg_0\equiv\0$. 
For each fixed $k\in\mathbb Z\backslash\{0\}$, a  solution to \eqref{3.21}--\eqref{3.22} is given by
\be
\bfu_k=\sum_{i=1}^3(\chi_{ki}\bfh_k^{(i)}+\sigma_{ki}\bfH_k^{(i)})\,,\ \ {\sf q}_k=\sum_{i=1}^3(\chi_{ki}p_k^{(i)}+\sigma_{ki}P_k^{(i)})\,,
\eeq{3.23}
with $(\bfh_k^{(i)},p_k^{(i)}), (\bfH_k^{(i)},P_k^{(i)})$ given in \lemmref{3.1}, and where $\bfchi_k, \bfsigma_k$ solve the equations
\be\ba{ll}\medskip
\i\,k\,{\bfchi_{k}}+\Sum{i=1}3\Int{\partial\Omega}{}\left[\chi_{ki}\mathbb T(\bfh_k^{(i)},{p}^{(i)}_k)+\sigma_{ki}\mathbb T(\bfH_k^{(i)},{P}^{(i)}_k)\right]\cdot\bfn=\bfcalf_k\,,
\\
\i\,k\,{\bfsigma_{k}}+\Sum{i=1}3\Int{\partial\Omega}{}\left[\chi_{ki}\bfx\times\mathbb T(\bfh_k^{(i)},{p}^{(i)}_k)+\sigma_{ki}\bfx\times\mathbb T(\bfH_k^{(i)},{P}^{(i)}_k)\right]\cdot\bfn=\bfcalg_k
\,.
\ea
\eeq{3.24}
Set
$$
\bfxi_k:=(\bfchi_{k},\bfsigma_{k})\in \mathbb C^6\,,\ \ \bfpzc F_k:=(\bfcalf_k,\bfcalg_k)\in \mathbb C^6\,.
$$
Then, with the notation of \lemmref{3.1},  \eqref{3.24} can be equivalently rewritten as
\be
\bfpzc B\cdot\bfxi_k=\bfpzc F_k\,,
\eeq{3.25}
where $\bfpzc B:=\i\,k\,\bfpzc I+\bfpzc A$, and $\bfpzc I$  is the $6\times 6$ identity matrix. By that lemma 
the matrix $\bfpzc B$ is invertible for all $k\in\mathbb Z$. Furthermore, from \eqref{3.A}, for all $\bfzeta\in\mathbb C^6$
we get
\be
\bfzeta^*\cdot\bfpzc B\cdot\bfzeta=\i\,k\,\left(|\bfzeta|^2+
\|\bfsf v\|_2^2\right)-\tau (\partial_1\bfsf v,\bfsf v^*)+2\|\mathbb D(\bfsf v)\|_2^2\,.
\eeq{3.26}
As a result, for any given $\bfpzc F_k$,  \eqref{3.25} has one and only one solution $\bfxi_k$.
We next dot-multiply both sides of \eqref{3.25} by $\bfxi_k^*$  and use \eqref{3.26} to deduce
$$
\i\,k\,\left(|\bfxi_k|^2+
\|\bfu_k\|_2^2\right)-\tau (\partial_1\bfu_k,\bfu_k^*)+2\|\mathbb D(\bfu_k)\|_2^2=(\bfpzc F_k,\bfxi_k^*)\,,
$$
which, by Cauchy--Schwarz inequality and \eqref{1.8} 
furnishes, in particular, the following estimates for all $k\in\mathbb Z\backslash\{0\}$
\be\ba{ll}\medskip 
\|\mathbb D(\bfu_k)\|_2^2\le \half\|\bfpzc F_k\|_2|\bfxi_k|\\ \medskip
|k|\,\|\bfu_k\|_2^2\le 2|\tau|^2\|\nabla\bfu_k\|_2^2 +2\|\bfpzc F_k\|_2|\bfxi_k|\le 4|\tau|^2\|\mathbb D(\bfu_k)\|_2^2+2\|\bfpzc F_k\|_2|\bfxi_k|\\
|k|\,|\bfxi_k|^2\le 2|\tau|\|\nabla\bfu_k\|_2\|\bfu_k\|_2 +\Frac{2}{|k|}\|\bfpzc F_k\|_2^2\le 2^{\frac32} |\tau|\|\mathbb D(\bfu_k)\|_2\|\bfu_k\|_2 +\Frac{2}{|k|}\,\|\bfpzc F_k\|_2^2
\,. 
\ea
\eeq{3.27}
Replacing \eqref{3.27}$_1$ into \eqref{3.27}$_2$, we obtain
\be
|k|\,\|\bfu_k\|_2^2\le c\,\|\bfpzc F_k\|_2|\bfxi_k|\,,
\eeq{3.28}
while using \eqref{3.27}$_1$ and \eqref{3.28} into \eqref{3.27}$_3$ along with Cauchy-Schwarz inequality implies 
\be
|k|^2\,|\bfxi_k|_2^2\le c\,\|\bfpzc F_k\|_2^2\,.
\eeq{3.29}
Combining \eqref{3.27}$_1$, \eqref{3.28}, \eqref{3.29} and \eqref{1.8}, and recalling \eqref{3.4}$_3$ and \eqref{3.23} we thus infer
\be
\sum_{|k|\ge1}\left[(|k|^2+1)\|\bfu_k\|_2^2+\|\nabla\bfu_k\|_2^2+\|D^2\bfu_k\|_{2}^2\right]\le c\sum_{|k|\ge1}(|k|^2+1)|\bfxi_k|^2\le c\,\|\bfpzc F\|_{L^2}^2\,.
\eeq{3.30}
Therefore, we may conclude that the quadruple $(\bfu,{\sf q},\bfchi_\sub{\bfu}\equiv\bfchi,\bfsigma_\sub{\bfu}\equiv\bfsigma)$ defined in \eqref{3.20} with $(\bfu_k,{\sf q}_k,\bfchi_k,\bfsigma_k)$ satisfying \eqref{3.21}--\eqref{3.22} is a solution to \eqref{3.17} in the class $\mathcal W_\sharp^2\times\calp^{1,2}\times W^{1,2}_\sharp\times W^{1,2}_\sharp$. Furthermore, \eqref{3.19} and \eqref{3.30} also entail the validity of the following inequality
$$
\|\bfu\|_{\mathcal W_\sharp^{2}}+\|{\sf q}\|_{\mathcal P^{1,2}}+\|\bfchi\|_{W^{1,2}}+\|\bfsigma\|_{W^{1,2}}\le c\,\Big(\|\bff\|_{\call_\sharp^2}+\|\bfF\|_{L^2}+\|\bfG\|_{L^2}\Big)\,.
$$
The existence proof is thus completed. 
The uniqueness property amounts to show that the problem
\be\ba{cc}\medskip\left.\ba{ll}\medskip
\partial_\s\bfw-\tau\,\partial_1\bfw=\Delta\bfw-\nabla {\sfp}\\ 
\Div\bfw=0\ea\right\}\ \ \mbox{in $\Omega\times[0,2\pi]$}\\ \medskip
\bfw|_{\partial\Omega}=\bfchi+\bfsigma\,;\\ \dot{\bfchi}+\Int{\partial\Omega}{}\mathbb T(\bfw,{\sf p})\cdot\bfn=\0\,,\ \ \dot{\bfsigma}+\Int{\partial\Omega}{}\bfx\times\mathbb T(\bfw,{\sf p})\cdot\bfn=\0
\ea
\eeq{3.31}
has only the zero solution in the specified function class. If we dot-multiply \eqref{3.31}$_1$ by $\bfw$, integrate by parts over $\Omega$ and use \eqref{3.31}$_3$, we get
$$
\half\ode{}t(\|\bfw(t)\|_2^2+|\bfchi(t)|^2+|\bfsigma(t)|^2)+2\|\mathbb D(\bfw(t))\|_2^2=0\,.
$$
Integrating both sides of this equation from $0$ to $2\pi$ and employing the $2\pi$-periodicity of the solution we easily obtain  $\|\mathbb D(\bfw(t))\|_2\equiv 0$ which, in turn, by the characterization of the space $\calv$ given in \lemmref{1.1}, immediately furnishes $\bfw\equiv\nabla\sfp\equiv\0$. The proof of the lemma is completed.\par\hfill$
\square$
\Br Concerning the generality of the domain $\Omega$, an observation similar
to that made in \remref{3.1} for \lemmref{3.1}, equally applies also to \lemmref{3.2}.
\ER{3.2}
\par
Let $\bfcalf\in\mathfrak H_\sharp$ where 
$$
\bfcalf=\left\{\ba{ll}\medskip\bff \ \mbox{in $\Omega$}\\
\bfF+\bfG\times\bfx\ \mbox{in $\Omega_0$}\ea\right.\,,
$$ 
and consider the operator equation
\be
\bfcalq_0(\bfw)=\bfcalf\,.
\eeq{chiappa}
By \lemmref{0.1}, \eqref{chiappa} is {\em equivalent} to the following problem (with $\bftau:=\lambda_\C\bftau_\C$, $\bfchi:=\hat{\bfw}_1$, $\bfsigma:=\hat{\bfw}_2$)
\be \ba{cc}\medskip\left.\ba{ll}\medskip
\zeta_0\,\partial_\s\bfw-\lambda_\C\bftau_\C\cdot\nabla\bfw-{\Delta}\bfw=\nabla \phi+\bff\\ 
\Div\bfw=0\ea\right\}\ \ \mbox{in $\Omega\times[0,2\pi]$}\,,\\ \medskip
\bfw=\bfchi+\bfsigma\times\bfx\ \ \mbox{at $\partial\Omega\times[0,2\pi]$}\,,\\
\,M\dot{\bfchi}+\Int{\partial\Omega}{}\mathbb T(\bfw,{\phi})\cdot\bfn-\bfF+   
\left(\mathcal I\,\dot{\bfsigma}+\Int{\partial\Omega}{}\bfy\times\mathbb T(\bfw,{\phi})\cdot\bfn-\bfG\right)\times\bfx=\0\, \ \ \mbox{in $\Omega_0\times[0,2\pi]$}\,.
\ea
\eeq{ciapa}
Since $\bfx$ is arbitrary in $\Omega_0$, we conclude that \eqref{ciapa}$_{4,5}$ are equivalent to \eqref{3.1}$_{4,5}$. Thus, 
in view of \lemmref{3.2}, we deduce the following important result
\Bl The operator $\bfcalq_0$ is a homeomorphism.
\EL{Hom}
This lemma allows us to prove the following theorem that represents the main result of this section.
\Bt
The operator $\bfcalq$ is Fredholm of index 0.
\ET{3.1}
{\em Proof.} We commence to notice that $\bfcalq=\bfcalq_0+\mathscr P\,\bfcalk$. Thus, by \lemmref{3.2},  the stated property will follow, provided we show that the map
$$
\bfcalc :\bfw\in\mathfrak W_\sharp^2\mapsto \bfcalk(\bfw)\in \mathfrak L^2_\sharp
$$
is compact. Let $\{\bfw_k\}$ be a bounded sequence in $\mathfrak W^{2}_\sharp$. This implies, in particular, that there is $M>0$ independent of $k$ such that ($\bfchi_k:=\widehat{(\bfw_k)}_1$)
\be
\|\bfw_k\|_{\calw_\sharp^2}+\|\bfchi_k\|_{W^{1,2}}\le M\,.
\eeq{gatto1}
We may then select sequences (again denoted by $\{\bfw_k,\bfchi_k\}$) and find $(\bfw_*,\bfchi_*)\in \mathcal W^{2}_\sharp\times W^{1,2}_\sharp$ such that
\be
\bfw_k\to{\bfw_*} \ \, \mbox{weakly in $\mathcal W^{2}_\sharp$\,;}\ \ \bfchi_k\to{\bfchi_*}\,  \  \mbox{strongly in $L^{\infty}(0,2\pi)$.}
\eeq{2.18}
Without loss of generality, we may take $\bfw_*\equiv\bfchi_*\equiv\0$.
We then have to show that
\be
\lim_{k\to\infty}\int_0^T\|\bfcalk(\bfw_k)\|_{2,\Omega}^2=0\,.
\eeq{gatto2}
From \eqref{2.18},  the compact embeddings $W^{2,2}(\Omega)\subset W^{1,4}(\Omega_R)\subset L^{2}(\Omega_R)$ for all $R>1$, and Lions-Aubin lemma we then have
\be
\int_{0}^{2\pi}\left(\|\bfw_k(t)\|_{2,\Omega_R}^2+\|\nabla\bfw_k(t)\|_{2,\Omega_R}^2\right){\rm d}t\to 0\ \ \mbox{as $k\to\infty$, for all $R>1$\,.}
\eeq{2.19}
Further, by \eqref{stoc},  
$$
\int_{0}^{2\pi} \|\bfv_\C\cdot\nabla\bfw_k(t)\|_{2,\Omega}^2
\le \|\bfv_\C\|_\infty^2\int_{0}^{2\pi} \|\nabla\bfw_k(t)\|_{2,\Omega_R}^2+ \|\bfv_\C\|_{4,\Omega^R}^2\int_{-\pi}^{\pi}\|\nabla\bfw_k(t)\|_{4,\Omega}^2\,,
$$
which, by \eqref{stoc}$_1$, \eqref{gatto1}, \eqref{2.19}  and the arbitrariness of $R$ furnishes
\be
\lim_{k\to\infty}\int_{0}^{2\pi} \|\bfv_\C\cdot\nabla\bfw_k(\tau)\|_{2}^2=0\,.
\eeq{2.21}
Likewise, 
$$
\Int{0}{2\pi} \|\bfw_k(t)\cdot\nabla\bfv_\C\|_{2,\Omega}^2
\le \|\nabla\bfv_\C\|_{\infty}^2\Int{0}{2\pi} \|\bfw_k(t)\|_{2,\Omega_R}^2 + 
\|\nabla\bfv_\C\|_{2,\Omega^R}^2\Int{0}{\pi}\|\bfw_k(t)\|_{2,\Omega}^2\,,
$$
so that, by \eqref{stoc}$_2$,  \eqref{2.18}$_1$, and \eqref{2.19} we deduce, as before,
\be
\lim_{k\to\infty}\int_{0}^{2\pi} \|\bfw_k(t)\cdot\nabla\bfv_\C\|_{2,\Omega}^2=0\,.
\eeq{2.22}
Finally,
$$
\int_{0}^{2\pi}\|\bfchi_k\cdot\nabla\bfv_\C\|_{2,\Omega}^2\le {2\pi}\, \|\bfchi_k\|_{L^\infty(0,2\pi)}^2\|\nabla\bfv_\C\|_2^2\,,
$$
which, by \eqref{stoc}$_2$ and \eqref{2.18}$_2$ furnishes
\be
\lim_{k\to\infty}\int_{0}^{2\pi}\|\bfchi_k\cdot\nabla\bfv_\C\|_{2,\Omega}^2=0
\eeq{3.38}
Combining \eqref{2.21}--\eqref{3.38} we thus arrive at \eqref{gatto2},
which completes the proof of the theorem.\par
\hfill$\square$
\setcounter{equation}{0}
\section{Sufficient Conditions for Time-Periodic Bifurcation}
The first objective of this section is to  rewrite \eqref{1.3} in an operator form of the type \eqref{Ar.5} and then, successively, employ \theoref{3.1_ar} to provide sufficient conditions for the occurence of time-periodic bifurcation for our problem.  Thus,
let
$$
\bfu(x,\s)=\bar{\bfu}(x)+\bfw(x,\s)\,,\ \ {\sf p}=\bar{\sf p}(x)+\phi(x,\s)\,,\ \ \bfgamma=\bar{\bfgamma}+\bfchi(\s)\,,\ \ \bfomega{(\s)}=\bar{\bfomega}+\bfsigma(\s)\,.
$$  
Then, \eqref{1.3} can be equivalently written in terms of the two sets of unknowns $(\bar{\bfu},\bar{\sf p},\bar{\bfgamma},\bar{\bfomega})$ and $(\bfw,\phi,\bfchi,\bfsigma)$ as follows
\be\ba{cc}\medskip\left.\ba{ll}\medskip-\lambda\bftau\cdot\nabla\bar{\bfu}+\lambda\,(\bfv\cdot\nabla\bar{\bfu}+(\bar{\bfu}-\bar{\bfgamma})\cdot\nabla\bfv)+\lambda\,(\bar{\bfu}-\bar{\bfgamma})\cdot\nabla\bar{\bfu}+\lambda\,\bar{(\bfw-\bfchi)\cdot\nabla\bfw}\\ \medskip\hspace*{6cm}=\Delta\bar{\bfu}-\nabla\bar{\sf p}\\
\Div\bar{\bfu}=0\ea\right\}\,  \mbox{in $\Omega$}\,,\\ \medskip
\bar{\bfu}=\bar{\bfgamma}+\bar{\bfomega}\times\bfx\ \ \mbox{at $\partial\Omega$}\,,\ \ \Lim{\mbox{\footnotesize $|\bfx|$}\to\infty}\bar{\bfu}(x)=\0\,,\\ 
\Int{\partial\Omega}{}\mathbb T(\bar{\bfu},\bar{\sf p})\cdot\bfn=\0\,,\ \   
\Int{\partial\Omega}{}\bfx\times\mathbb T(\bar{\bfu},\bar{\sf p})\cdot\bfn=\0\, \,,
\ea
\eeq{1.4}
and
\be\ba{cc}\medskip\left.\ba{ll}\medskip\zeta\,\partial_\s\bfw-\lambda\bftau\cdot\nabla\bfw+\lambda(\bfv\cdot\nabla\bfw+(\bfw-\bfchi)\cdot\nabla\bfv)+\lambda\big[(\bfw-\bfchi)\cdot\nabla\bar{\bfu}\\ \medskip\hspace*{2cm}+(\bar{\bfu}-\bar{\bfgamma})\cdot\nabla\bfw+(\bfw-\bfchi)\cdot\nabla\bfw-\bar{(\bfw-\bfchi)\cdot\nabla\bfw}\big]=\Delta\bfw-\nabla{\phi}\\
\Div\bfw=0\ea\right\}\,  \mbox{in $\Omega\times [0,2\pi]$}\,,\\ \medskip
\bfw=\bfchi+\bfsigma\times\bfx\ \ \mbox{at $\partial\Omega\times [0,2\pi]$}\,,\ \ \Lim{\mbox{\footnotesize $|\bfx|$}\to\infty}\bfw(x,\s)=\0\,,\ \ \s\in[0,2\pi]\,,\\ 
\zeta\,M\dot{\bfchi}+\Int{\partial\Omega}{}\mathbb T(\bfw,{\phi})\cdot\bfn=\0\,,\ \   
\zeta\,\mathcal I\,\dot{\bfsigma}+\Int{\partial\Omega}{}\bfx\times\mathbb T(\bfw,{\phi})\cdot\bfn=\0\, \ \ \mbox{in $[0,2\pi]$}\,.
\ea
\eeq{1.5}
where, for simplicity, we have suppressed the subscript \mbox{{\footnotesize 0}}.
Let $\bfv_\C$ be the weak solution to \eqref{1.2} at $\lambda=\lambda_\C$, and let $\bftau_\C$ be the associated translational velocity. We make the assumption that both $\bftau$ and $\bar{\bfgamma}$ are directed along the direction $\bfe_\C:=\bftau_\C/|\bftau_\C|$, and write $\bftau=\tau\,\bfe_\C$, $\bar{\bfgamma}=\bar{\gamma}\,\bfe_\C$. As a result, by \theoref{2.2} we know that, under the hypothesis \eqref{2.13},  at $\lambda=\lambda_\C$ there exists an analytic family of weak solutions $\bfv=\bfv(\lambda)$ such that $\bfv(\lambda_\C)=\bfv_\C$.  Thus, setting $\mu:=\lambda-\lambda_\C$, \eqref{1.4} can be rewritten as follows
\be\ba{cc}\medskip\left.\ba{ll}\medskip
-\Delta\bar{\bfu}-\lambda_\C\tau_\C\bfe_\C\cdot\nabla\bar{\bfu}+\lambda_\C\,(\bfv_\C\cdot\nabla\bar{\bfu}+(\bar{\bfu}-\bar{\gamma}\,\bfe_\C)\cdot\nabla\bfv_\C)=\bfN_1(\bar{\bfu},\bfw,\mu)-\nabla\bar{\sf p}\\
\Div\bar{\bfu}=0\ea\right\}\,  \mbox{in $\Omega$}\,,\\ \medskip
\bar{\bfu}=\bar{\gamma}\,\bfe_\C+\bar{\bfomega}\times\bfx\ \ \mbox{at $\partial\Omega$}\,,\ \ \Lim{\mbox{\footnotesize $|\bfx|$}\to\infty}\bar{\bfu}(x)=\0\,,\\ 
\Int{\partial\Omega}{}\mathbb T(\bar{\bfu},\bar{\sf p})\cdot\bfn=\0\,,\ \   
\Int{\partial\Omega}{}\bfx\times\mathbb T(\bar{\bfu},\bar{\sf p})\cdot\bfn=\0\, \,,
\ea
\eeq{7.3}
where
\be
\ba{rl}\medskip
\bfN_1(\bar{\bfu},\bfw,\mu):=&\!\!\!\!-[(\tilde{\mu}\,\tilde{\tau}(\mu)-\lambda_\C\tau_\C)\,\bfe_\C+\tilde{\mu}\,\tilde{\bfv}(\mu)-\lambda_c\bfv_\C]\cdot\nabla\bar{\bfu}\\ \medskip
&\!\!\!\!
-\lambda_\C(\bar{\bfu}-\bar{\gamma}\bfe_\C)\cdot\nabla(\tilde{\bfv}(\mu)-\bfv_\C)-\mu(\bar{\bfu}-\bar{\gamma}\bfe_\C)\cdot\nabla\tilde{\bfv}(\mu)\\
&\!\!\!\!-\tilde{\mu}[\,(\bar{\bfu}-\bar{\gamma}\,\bfe_\C)\cdot\nabla\bar{\bfu}+\bar{(\bfw-\bfchi)\cdot\nabla\bfw}]\,,
\ea
\eeq{7.4}
and $\tilde{f}(\mu):=f(\mu+\lambda_\C)$. Likewise, \eqref{1.5} can be rewritten as follows
\be\ba{cc}\medskip\left.\ba{ll}\medskip\zeta\,\partial_\s\bfw-\lambda_\C\tau_\C\bfe_\C\cdot\nabla\bfw+\lambda_\C(\bfv_\C\cdot\nabla\bfw+(\bfw-\bfchi)\cdot\nabla\bfv_\C)-\Delta\bfw\\
\hspace*{8cm}=\bfN_2(\bar{\bfu},\bfw,\mu)-\nabla{\phi}\\
\Div\bfw=0\ea\right\}\,  \mbox{in $\Omega\times [0,2\pi]$}\,,\\ \medskip
\bfw=\bfchi+\bfsigma\times\bfx\ \ \mbox{at $\partial\Omega\times [0,2\pi]$}\,,\ \ \Lim{\mbox{\footnotesize $|\bfx|$}\to\infty}\bfw(x,\s)=\0\,,\ \ \s\in[0,2\pi]\,,\\ 
\zeta\,M\dot{\bfchi}+\Int{\partial\Omega}{}\mathbb T(\bfw,{\phi})\cdot\bfn=\0\,,\ \   
\zeta\,\mathcal I\,\dot{\bfsigma}+\Int{\partial\Omega}{}\bfx\times\mathbb T(\bfw,{\phi})\cdot\bfn=\0\, \ \ \mbox{in $[0,2\pi]$}\,,
\ea
\eeq{7.5}
where
\be
\ba{rl}\medskip
\bfN_2(\bar{\bfu},\bfw,\mu):=&\!\!\!\!-[(\tilde{\mu}\,\tilde{\tau}(\mu)-\lambda_\C\tau_\C)\,\bfe_\C+\tilde{\mu}\,\tilde{\bfv}(\mu)-\lambda_c\bfv_\C]\cdot\nabla{\bfw}\\ \medskip
&\!\!\!\!
-\lambda_\C({\bfw}-{\bfchi})\cdot\nabla(\tilde{\bfv}(\mu)-\bfv_\C)-\mu({\bfw}-{\bfchi})\cdot\nabla\tilde{\bfv}(\mu)\\
&\!\!\!\!-\tilde{\mu}\,\big[(\bfw-\bfchi)\cdot\nabla\bar{\bfu}+(\bar{\bfu}-\bar{\bfgamma})\cdot\nabla\bfw+(\bfw-\bfchi)\cdot\nabla\bfw-\bar{(\bfw-\bfchi)\cdot\nabla\bfw}\big]\,,
\ea
\eeq{7.6}
Let
$$
\tilde{\bfN}_2=\left\{\ba{cc}\medskip {\bfN}_2\,\ \mbox{in $\Omega$}\\ \0\,\ \,\mbox{in $\Omega_0$}\ea\right.\,
.
$$
The following result holds.
\Bl The operators
\be\ba{ll}\medskip
\bfscr N_1:(\bar{\bfu},\bfw,\mu)\in {\mathscr X}_\C(\real^3)\times\mathfrak W_\sharp^2\times\real\mapsto \bfN_1(\bar{\bfu},\bfw,\mu)\in\mathcal V^{-1}(\real^3)
\\
\tilde{\bfscr N}_2:(\bar{\bfu},\bfw,\mu)\in {\mathscr X}_\C(\real^3)\times\mathfrak W_\sharp^2\times\real\mapsto \tilde{\bfN}_2(\bar{\bfu},\bfw,\mu)\in\mathfrak L^2_\sharp
\ea
\eeq{7.7}
are well defined.
\EL{7.1}
{\em Proof.} Since $\bar{\bfu},\bfv,\bfv_\C\in\mathscr X_\C$, it follows at once that $\bfe_\C\cdot\nabla\bar{\bfu}, \bfe_\C\cdot\nabla{\bfv}, \bfe_\C\cdot\nabla{\bfv_\C}\in \mathcal V^{-1}$. Moreover, by \lemmref{1.2} we also have $\bar{\bfu},\bfv,\bfv_\C\in L^4(\real^3)$, which, by using integration by parts, implies $\bfv\cdot\nabla\bar{\bfu},\bfv_\C\cdot\nabla\bar{\bfu},\bar{\bfu}\cdot\nabla{\bfv},\bar{\bfu}\cdot\nabla{\bfv_\C},\bar{\bfu}\cdot\nabla\bar{\bfu}\in\calv^{-1}(\real^3)$ as well. Finally, observing that $\bfw\in L^2(L^2\cap L^4)$ and $\bfchi\in L^\infty(0,2\pi)$ we easily show that $\bar{(\bfw-\bfchi)\cdot\nabla\bfw}\in \calv^{-1}$, which concludes the proof of \eqref{7.7}$_1$. By known  embedding theorems \cite[Theorem 2.1]{Solo} it follows that $\bfw\in L^2(D^{1,4})$. Thus, the validity of \eqref{7.7}$_2$ can be established along the same lines used to show \eqref{7.7}$_1$. We will omit the details. The proof of the lemma is completed.\par\hfill$\square$\par
In view of \lemmref{7.1}, and  \eqref{L1}, \eqref{L1_1}, and \eqref{Q}, it follows at once that, setting $\bfscr N_2=\mathscr P\,\tilde{\bfscr N}_2$, the coupled problem \eqref{7.3}--\eqref{7.6} can be written as operator equations:
\be\ba{ll}\medskip 
\bfscr L_1(\bar{\bfu})=\bfscr N_1(\bar{\bfu},\bfw,\mu)\ \ \mbox{in $\calv^{-1}(\real^3)$}\,,\\
\zeta\,\partial_\s\bfw+\bfscr L_2(\bfw)=\bfscr N_2(\bar{\bfu},\bfw,\mu)\ \ \mbox{in $\mathfrak L^2_\sharp$}\,
\ea
\eeq{7.8}
which are  in the form \eqref{Ar.5}. We shall next check  how the assumptions of \theoref{3.1_ar} can be satisfied in our case. We begin to notice that \theoref{3.1} secures  (H3). In addition, both hypotheses (H1) and (H4) are verified if we assume
\tag{$\calh1$}
\be
{\sf N}\,[\bfscr L_1]=\{\0\}\,.
\eeq{h1}
In fact, by  \theoref{2.2}, $\bfscr L_1$ is Fredholm of index 0, so that \eqref{h1} implies (H1). Moreover, if \eqref{h1} holds, then  --taking into account that the nonlinear operators $\bfscr N_i$, $i=1,2$,  are (at most) quadratic in $(\bar{\bfu},\bfw$)-- again by \theoref{2.2}, we deduce the validity  of (H4). Next,  we assume
\tag{$\calh2$}
\be
\mbox{$\nu_0:=\i\,\zeta_0$, $\zeta_0\neq0$,  is a simple eigenvalue of $\bfscr L_2$, and  $k\,\nu_0\in {\sf P}(\bfscr L_2)$, for all $k\in\nat\backslash\{0,1\}$\,. 
}
\eeq{h2}
This, in view of  \theoref{5.1}, guarantees  assumption (H2). By a straightforward calculation, we can sho that, in our case, the operator $S_{011}$ is given by
$$
(\tau_c\,\bfe_\C+\bfv_\C)\cdot\nabla\bfw+\lambda_\C(\tau^\prime(\lambda_\C)\bfe_c+\bfv^\prime(\lambda_\C))\cdot\nabla\bfw+(\bfw-\bfchi)\cdot\nabla(\bfv_\C+\bfv^\prime(\lambda_\C))\,,
$$
where the prime means differentiation with respect to $\mu$. So, denoting by $\nu=\nu(\mu)$ the eigenvalues of $\bfscr L_2+\mu\,S_{011}$, by  \cite[Proposition 79.15 and Corollary 79.16]{Z1}) we have that in a neighborhood of $\mu=0$ the map $\mu\mapsto\nu(\mu)$ is well defined and of class $C^\infty$. This justifies our last assumption:
\tag{$\calh3$}
\be
\Re[\nu^{\prime}(0)]\neq 0\,.
\eeq{h3}
\par We are now in a position to employ \theoref{3.1_ar} to obtain the following main result.
\renewcommand{\theequation}{\arabic{section}.\arabic{equation}}\setcounter{equation}{8}
\Bt
Suppose  \eqref{h1}--\eqref{h3} hold. Let $\bfw_0$ be the normalized eigenvector of $\bfscr L_2$ corresponding to the eigenvalue $\nu_0$, and set
$
\bfw_1:=\Re[\bfw_0\,{\rm e}^{-{\rm i}\,\s}].$    
Then, the following properties are valid. \smallskip\\
{\rm (a)} {\rm Existence.} There are analytic families
\be
\big(\bar{\bfu}(\varepsilon),\bfw(\varepsilon),\zeta(\varepsilon),\mu(\varepsilon)\big)\in \mathscr X_\C\times \mathfrak W^{2}_\sharp\times \real_+\times\real
\eeq{fam1}
satisfying \eqref{7.8}, for all real $\varepsilon$ in a neighborhood $\mathcal I(0)$ of 0, and such that
\be
\big(\bar{\bfu}(\varepsilon),\bfw(\varepsilon)-\varepsilon\,\bfw_1,\zeta(\varepsilon),\mu(\varepsilon)\big)\to (0,0,\zeta_0,0)\ \ \mbox{as $\varepsilon\to 0$}\,.
\eeq{Ar.101}
\par\noindent
{\rm (a)} {\rm Uniqueness.}
There is a neighborhood  $$U(0,0,\zeta_0,0)\subset \mathscr X_\C\times \mathfrak W^{2}_\sharp\times \real_+\times \real$$ such that every (nontrivial) $2\pi$-periodic solution to \eqref{7.8},  lying in $U$ must coincide, up to a phase shift, with a member of the family \eqref{fam1}.
\smallskip\par\noindent
{\rm (a)} {\rm Parity.}  The functions $\zeta(\varepsilon)$ and $\mu(\varepsilon)$ are even:
$$
\zeta(\varepsilon)=\zeta(-\varepsilon)\,,\ \ \mu(\varepsilon)=\mu(-\varepsilon)\,,\ \ \mbox{for all $\varepsilon\in\cali(0)$\,.} 
$$
Consequently, the bifurcation due to these solutions is either subcritical or supercritical, a two-sided bifurcation being excluded.\footnote{Unless $\mu\equiv 0$.}
\ET{7.1}
\setcounter{equation}{0}
\section{On the Motion of the Sphere in the Time-Periodic Regime}
As we mentioned in the introductory section, experimental and numerical tests show that, in the transition from steady to time-periodic motion, the trajectory of the center of mass, $G$, of the sphere changes from a rectilinear to a zigzag path. Objective of this section is to study in more details the motion of the sphere in the time-periodic regime and, in particular, to derive necessary and sufficient conditions for the occurrence of this sideway oscillatory behavior.  
\par
To this end, we begin to observe that, according to \theoref{7.1}, in the neighborhood of $\lambda_c$, namely, for sufficiently small $\varepsilon$,  the oscillatory part of the solution, $\bfw$, behaves like the corresponding solution to the linear problem $\bfcalq(\bfw_1)=\0$, that is, $\bfw_1:=\Re[\bfw_0\,{\rm e}^{-{\rm i}\,\s}]$, with $\bfw_0$ eigenvector of $\bfscr L_2$ corresponding to the eigenvalue $-\i\,\zeta_0$. Therefore, in such a neighborhood, the oscillatory component of the velocity of $G$ will have the same kinematic characteristics of the translational velocity, $\bfchi_0$, and angular velocity, $\bfsigma_0$ associated to $\bfw_0$. We now recall that the equation $\bfscr L_2(\bfw_0)+\i\,\zeta_0\,\bfw_0=\0$ is equivalent to the following set of equations
\be\ba{cc}\medskip\left.\ba{ll}\medskip-\Delta\bfw_0-\lambda_\C\bftau_\C\cdot\nabla\bfw_0+\lambda_\C\,L(\bfw_0)+\nabla{\mathfrak p}+\i\,\zeta\,\bfw_0=\0\\
\Div\bfw_0=0\ea\right\}\,  \mbox{in $\Omega$}\,,\\ \medskip
\bfw_0=\bfchi_0+\bfsigma_0\times\bfx\ \ \mbox{at $\partial\Omega$}\,,\\ 
\i\,\zeta\,M\,{\bfchi_0}+\Int{\partial\Omega}{}\mathbb T(\bfw_0,{\mathfrak p})\cdot\bfn=\0\,,\ \   
\i\,\zeta\,\mathcal I\,{\bfsigma_0}+\Int{\partial\Omega}{}\bfx\times\mathbb T(\bfw_0,{\mathfrak p})\cdot\bfn=\0\,,
\ea
\eeq{8.1}
where  
\be
L(\bfw_0):= \bfv_\C\cdot\nabla\bfw_0+(\bfw_0-\bfchi_0)\cdot\nabla\bfv_\C\,.
\eeq{8.2}
Thus, assuming the gravity directed along $\bfe_1$, an oscillatory motion of $G$ in the neighborhood of $\lambda=\lambda_\C$ will take place if and only if $(\chi_0)_2\bfe_2+(\chi_0)_3\bfe_3\neq\0$. In the remaining part of this section we shall furnish a characterization of the expression of $\bfchi_0$ and $\bfsigma_0$ that, in particular, will provide the desired property.  
\par
Let us introduce the pairs $(\bfh^{(i)},p^{(i)})$ and $(\bfH^{(i)},P^{(i)})$ in $W^{2,2}(\Omega)\times D^{1,2}(\Omega)$, solutions to the following problems ($i=1,2,3$): 
\be
\ba{cc}\medskip\left.\ba{ll}\medskip
-{\rm i}\,\zeta_0\,\bfh^{(i)}+\lambda_\C\bftau_\C\cdot\nabla\bfh^{(i)}=\Div\mathbb T(\bfh^{(i)},p^{(i)})\\
\Div \bfh^{(i)}=0\ea\right\}\ \ \mbox{in $\Omega$}\\
\bfh^{(i)}=\bfe_i\,\ \ \mbox{at $\partial\Omega$}\,,
\ea
\eeq{8.3}
and
\be
\ba{cc}\medskip\left.\ba{ll}\medskip
-{\rm i}\,\zeta_0\,\bfH^{(i)}+\lambda_\C\bftau_\C\cdot\nabla\bfH^{(i)}=\Div \mathbb T(\bfH^{(i)},P^{(i)})\\
\Div \bfH^{(i)}=0\ea\right\}\ \ \mbox{in $\Omega$}\\
\bfH^{(i)}=\bfe_i\times\bfx\ \ \mbox{at $\partial\Omega$}\,.
\ea
\eeq{8.4}
Moreover, consider the matrices $\hat{\mathbb K}, \hat{\mathbb A}, \hat{\mathbb P}$, and $\hat{\mathbb S}$, defined by ($i,j=1,2,3$)
\be\ba{cc}\ms
(\hat{\mathbb K})_{ij}=\IdS(\mathbb T(\bfh^{(i)*},p^{(i)*})\cdot\bfn)_j,\\ \ms
(\hat{\mathbb A})_{ij}=\IdS(\bfx\times\mathbb T(\bfH^{(i)*},P^{(i)*})\cdot\bfn)_j\\ \ms
(\hat{\mathbb P})_{ij}=\IdS(\bfx\times\mathbb T(\bfh^{(i)*},p^{(i)*})\cdot\bfn)_j\\ (\hat{\mathbb S})_{ij}=\IdS(\mathbb T(\bfH^{(i)*},P^{(i)*})\cdot\bfn)_j\,,
\ea
\eeq{8.5}
where, we recall, ${}^*$ means complex conjugate.
The existence of the above pairs in the specified function class is guaranteed by \lemmref{3.1}. Furthermore, again from this lemma, we know that the matrices $\hat{\mathbb K}+\i\,\lambda \mathds{1}$ and $\hat{\mathbb A}+\i\,\mu \mathds{1}$ are invertible, for all $\lambda,\mu\in\real$, as well as the block matrix
$$
\hat{\bfpzc A}=\left(\ba{cc}\hat{\mathbb K}+\i\,\lambda \mathds{1} & \hat{\mathbb P}\\
\hat{\mathbb S} & \hat{\mathbb A}+\i\,\mu \mathds{1}\ea\right)
\,.
$$ 
If we dot-multiply both sides of \eqref{8.1}$_1$ by $\bfh^{(i)*}$,  integrate by parts over $\Omega$, and employing \eqref{8.1}$_{2,3,4}$, we get
\be
\i\,\zeta_0(\bfw_0,\bfh^{(i)*})-\lambda_\C(\bftau_\C\cdot\nabla\bfw_0,\bfh^{(i)*})+2[\mathbb D(\bfw_0),\mathbb D(\bfh^{(i)*})]=-\i\,M\,\zeta_0\bfchi_0\cdot\bfe_i-\lambda_\C(L(\bfw_0),\bfh^{(i)*})\,.
\eeq{8.6}
Similarly, taking first the complex conjugate of \eqref{8.3}$_1$, then dot-multiplying it by $\bfw_0$,  integrating by parts over $\Omega_0$, and using \eqref{8.5}, we deduce
\be
\i\,\zeta_0(\bfw_0,\bfh^{(i)*})-\lambda_\C(\bftau_\C\cdot\nabla\bfw_0,\bfh^{(i)*})+2[\mathbb D(\bfw_0),\mathbb D(\bfh^{(i)*})]=[\hat{\mathbb K}\cdot\bfchi_0+\hat{\mathbb P}\cdot\bfsigma_0]_i \,.
\eeq{8.7}
From \eqref{8.6} and \eqref{8.7} we conclude (summation over repeated indeces)
\be
\tilde{\mathbb K}\cdot\bfchi_0+\hat{\mathbb P}\cdot\bfsigma_0=-\lambda_\C(L(\bfw_0),\bfh^{(i)*})\,\bfe_i:=\bfpzc F\,,
\eeq{8.8}
where $\tilde{\mathbb K}:=\i\,M\,\zeta_0\,\mathds{1}+\hat{\mathbb K}$. 
Likewise, we can show that
\be
\tilde{\mathbb A}\cdot\bfsigma_0+\hat{\mathbb S}\cdot\bfchi_0=-\lambda_\C(L(\bfw_0),\bfH^{(i)*})\,\bfe_i:=\bfpzc G\,,
\eeq{8.9}
with $\tilde{\mathbb A}:=\i\,\cali\,\zeta_0\,\mathds{1}+\hat{\mathbb A}$. From \eqref{8.8} and \eqref{8.9} we infer
\be 
\bfchi_0=\mathbb H\cdot(\bfpzc F+\tilde{\mathbb K}\cdot\tilde{\mathbb A}^{-1}\cdot\bfpzc G)\,,\ \ \bfsigma_0=\mathbb M\cdot(\bfpzc G+\tilde{\mathbb A}\cdot\tilde{\mathbb K}^{-1}\cdot\bfpzc F)\,,
\eeq{8.10}
where
$$
\mathbb H:= (\tilde{\mathbb K}-\hat{\mathbb P}\cdot\tilde{\mathbb A}^{-1}\cdot\hat{\mathbb S})^{-1}\,,\ \ \mathbb M:= (\tilde{\mathbb A}-\hat{\mathbb S}\cdot\tilde{\mathbb K}^{-1}\cdot\hat{\mathbb P})^{-1}   \,.
$$
Notice that both $\mathbb H$ and $\mathbb M$ exist, because $\hat{\bfpzc A}$, $\tilde{\mathbb K}$, and $\tilde{\mathbb A}$ are invertible.
\par
From \eqref{8.10} we can then derive the following result.
\Bt Suppose the assumptions of \theoref{7.1} hold. Then, an oscillatory motion of the center of mass $G$ of $\mathscr S$  in the neighborhood of $\lambda=\lambda_\C$ may occur if and only if 
$$
[\mathbb H\cdot(\bfpzc F+\tilde{\mathbb K}\cdot\tilde{\mathbb A}^{-1}\cdot\bfpzc G)]_2\bfe_2+
[\mathbb H\cdot(\bfpzc F+\tilde{\mathbb K}\cdot\tilde{\mathbb A}^{-1}\cdot\bfpzc G)]_3\bfe_3\neq \0\,.
$$

\ET{8.1}

\ed

Notice that any solution in the class $\mathfrak B$ to \eqref{4.1}--\eqref{4.2} with $\bfxi_\sub{\bfu}=\0$ is necessarily a solution to \eqref{5.4}--\eqref{5.5} in the class $\mathfrak C_0$. We shall show, however, that for ``almost all" $\calf\in \mathfrak L$ problem \eqref{5.4}--\eqref{5.5} is {\em not} well-posed in the class $\mathfrak C_0$, and this will imply that for ``almost all" $\calf\in \mathfrak L$ of ``small" norm, the corresponding solution to \eqref{4.1}--\eqref{4.2}  must have $\bfxi_\sub{\bfu}\neq\0$. In turn, the above lack of well-posedness  is secured if we show that the operator $\calm$ is Fredholm of {\em negative} index. Since from \lemmref{Op} and \lemmref{1.6} it follows that the index of $\cala$, ${\rm ind}\,\cala$, is $-3$. in order to show the same property for $\calm$ it suffices to prove that, at any given $\bfU_0\in \mathfrak C_0$,  the Fr\'echet derivative of $\caln+\call$ is compact.

Before performing this study, however, we would like to emphasize that whatever the shape of $\mathscr B$, and physical properties of $\mathscr B$ and $\mathscr L$, we can always find a boundary distribution that {\em does not} produce a non-zero net motion. This implies that, for a given body $\mathscr B$, some restrictions on the data are indeed necessary to guarantee  self-propulsion.  
Actually, let $\psi=\psi(x)$,  $x\in\Omega$, be a scalar function satisfying the following Neumann problem
\be
\Delta\psi=0\ \ \mbox{in $\Omega$}\,,\ \ \pde{\psi}{n}=f\ \ \mbox{at $\partial\Omega$}\ \ \ \mbox{with}\ \,\int_{\partial\Omega} f=0\,.
\eeq{p5}
Under condition \eqref{p5}$_3$ on the data, it is then well known that (for example, \cite[Exercise V.3.6]{Gab})
\be
D^\alpha \psi(x)=O(|x|^{-2-|\alpha|})\ \ \mbox{as $|x|\to\infty$}\,,\ \ |\alpha|\ge 0.
\eeq{p6}
Next, let $a=a(t)$ be a smooth $T$-periodic function with $\bar{a}=0$, and set
\be
\bfv(x,t):=a(t)\,\nabla\psi(x)\,.
\eeq{nb}
It is immediately checked that, in view of \eqref{p5}$_1$ and \eqref{p6}, the field $\bfv$ thus defined satisfies \eqref{SE}$_{1,2,4}$ with $ {\bfscr H}\equiv\0$, provided we choose 
\be 
p= -\dot{a}(t)\psi+a(t)\bfgamma(t)\cdot\nabla\psi-\half a^2(t)|\nabla\psi|^2\,.
\eeq{nb1}
Let us now turn to the surface integral in \eqref{SE}$_5$ (with $\bfscr{\bfscr H}\equiv\0$) that, with $\bfv$ and $p$ defined above becomes
\be\ba{ll}\smallskip
\Int{\partial\Omega}{}\Big[2a\,\nu\,\nabla(\nabla\psi)\cdot\bfn-\dot{a}\psi\,\bfn-a\,(\bfgamma\cdot\nabla\psi\,\bfn-\nabla\psi\,\bfgamma\cdot\bfn)\\
\hspace*{2.7cm}+a^2(\half|\nabla\psi|^2\bfn-\nabla\psi\nabla\psi\cdot\bfn)\Big]
:=2a\,\nu\,\bfI_1-\dot{a} \bfI_2+a\,\bfI_3+a^2\bfI_4\,.
\ea
\eeq{p7}
By a simple integration by parts and \eqref{p5}$_1$, we find, for $i=1,2,3$,
$$ 
I_{1i}=\int_{\partial\Omega}\partial_i\partial_j\psi\,n_j=\int_{\Omega_R}\partial_i(\partial_j\partial_j\psi)-\int_{\partial B_R}\partial_i\partial_j\psi\,n_j=-\int_{\partial B_R}\partial_i\partial_j\psi\,n_j\,,
$$
so that, letting $R\to\infty$ into this relation and using \eqref{p6}, we find
\be
\bfI_1=\0\,.
\eeq{p8}
We next observe that
$$
\bfgamma\cdot\nabla\psi\,\bfn-\nabla\psi\,\bfgamma\cdot\bfn=(\bfn\times\nabla\psi)\times \bfgamma\,,
$$
and that
$$
\int_{\partial\Omega}\bfn\times\nabla\psi=-\int_{\partial B_R}\bfn\times\nabla\psi+\int_{\Omega_R}\nabla\times(\nabla\psi)=-\int_{\partial B_R}\bfn\times\nabla\psi\,.
$$
As a result, letting $R\to\infty$ in the latter and employing \eqref{p6} we deduce
\be
\bfI_3=\0\,.
\eeq{p9}
Finally, again by integration by parts, $i=1,2,3$,
$$\ba{rl}\medskip
I_{4i}=&\!\!\!\!\Int{\Omega_R}{}\big[\half\partial_i|\nabla\psi|^2-\partial_\ell(\partial_i\psi\,\partial_\ell\psi)\big]-\Int{\partial B_R}{}\big(\half|\nabla\psi|^2n_i-\partial_i\psi\partial_\ell\psi\,n_\ell\big)\\
=&\!\!\!\!
-\Int{\partial B_R}{}\big(\half|\nabla\psi|^2n_i-\partial_i\psi\partial_\ell\psi\,n_\ell\big)
\ea
$$
and again by \eqref{p6}, we may let $R\to\infty$ in this relation and infer
\be
\bfI_4=\0\,.
\eeq{p10}
Collecting \eqref{p7}--\eqref{p10}, we then conclude that  a solution to \eqref{SE}  (with $\bfscr H\equiv\0$) is given by the fields $\bfv,p$ defined in \eqref{nb}--\eqref{nb1} together with 
$$
\bfgamma(t):=\frac{a(t)}M\int_{\partial\Omega}\psi\, \bfn\,,\ \ \bfv_*:=a(t)\,\nabla\psi-\bfgamma\,.
$$
Since, of course, $\bar{\bfgamma}=0$, no net motion occurs. Notice that $(\bfv,p,\bfgamma)$ is in the functional class $\mathcal W^{2,\frac32}_\sharp\times\calp^{1,\frac32}\times W^{1,2}_\sharp(0,T)$. 
 

Set 

Then, with the help of \lemmref{7.1}, we shall  show that the functions
\be
\bfv:=\delta\,\bfV_0\,,\ \ p:=a(t)\bfgamma(t)\cdot\nabla\psi-\half\,a^2(t)|\nabla\psi|^2+\delta\,P_0\,,\ \ \bfgamma:=\delta\bfchi_0
\eeq{VV}
satisfy \eqref{SE} with  ${\bfscr H}\equiv\0$ and $\bfv_*:=\delta\,\bfV_*\equiv \bfsf V_*$ defined in \lemmref{7.1}.

to the whole of $\real^3$ by a procedure entirely analogous to that used earlier on. Precisely, let $\bfscr S:=\nabla\bfP$, with $\bfP$ solving \eqref{NP} with $\bfscr F\equiv\bfscr H$ and set
$$
\tilde{\bfU}=\left\{\ba{ll}\medskip\bfU(x)\ &\mbox{if $x\in\Omega$}\\ 
\bfzeta\ &\mbox{if $x\in\Omega_0$}\ea\right.\,,\  \ \tilde{P}=\left\{\ba{ll}\medskip {P}(x)\ &\mbox{if $x\in\Omega$}\\ 
0\ &\mbox{if $x\in\Omega_0$}\ea\right.\,,\ \ \tilde{{\bfu}}(x)=\left\{\ba{ll}\medskip\bfu(x)\ &\mbox{if $x\in\Omega$}\\ 
{\bfxi_\sub{\bfu}}\ &\mbox{if $x\in\Omega_0$}\ea\right.\,,
$$ 
and
$$
\tilde{{\bfscr H}}(x)=\left\{\ba{ll}\medskip\bfscr{H}(x)\ &\mbox{if $x\in\Omega$}\\ 
{\bfscr S}(x)\ &\mbox{if $x\in\Omega_0$}\ea\right.\,.
$$
Then, arguing exactly as we did before, we show that $(\tilde{\bfU},\tilde{P})$ is a distributional solution to the problem
\be\left.\ba{ll}\medskip
\Delta\tilde{\bfU}+\bfxi_\sub{\bfu_1}\cdot\nabla\tilde{\bfU}=-\bfzeta\cdot\nabla\tilde{\bfu}+\nabla \tilde{P}+\Div\tilde{\bfscr H}\\ \Div\bfU=0
\ea\right\}\ \ \mbox{in $\real^3$}\,.
\eeq{St00}
\Bl
Assume   
\be
\Div  \bfscr{\bfscr G}\in \sfL^2\,, \ \ \bfb\in \call_\sharp^{2,\frac32}\,,\ \ \bfv_*\in W^{1,2}_\sharp\,,\ \ \bfg\in L^2_\sharp(0,T)\,.
\eeq{4.10}
Then, the following inequalities hold with $\bfU:=(\bfphi,\bfpsi,\bfchi_\sub{\bfpsi})$, 
\be\ba{ll}\medskip
\|\Div\bfscr{\bfscr F}(\bfU)\|_{\sfL^2}\le c\Big(\|\bfU\|_{\mathfrak B}^2+\|\bfU\|_{\mathfrak B}\|\bfv_*\|_{W^{1,2}(W^{\frac32,2}(\partial\Omega))}+\|\bfv_*\|_{W^{1,2}(W^{\frac32,2}(\partial\Omega))}^2+\|\Div\bfscr{\bfscr G}\|_{\sfL^2}\Big)\\ \medskip
\|\bff(\bfU)\|_{\call^{2,\frac32}}\le c\Big(\|\bfU\|_{\mathfrak B}^2\!\!+\!\|\bfv_*\|_{W^{1,2}(W^{\frac32,2}(\partial\Omega))}(1+\|\bfU\|_{\mathfrak B})+\|\bfv_*\|_{W^{1,2}(W^{\frac32,2}(\partial\Omega))}^2\!+\!\|\bfb\|_{\call^{2,\frac32}}\Big)\\
\|\bfF(\bfchi_\sub{\bfpsi},\bfxi_\sub{\bfphi})\|_{L^2(0,T)}\le
c\Big(\|\bfv_*\|_{W^{1,2}(W^{\frac32,2}(\partial\Omega))}(\|\bfU\|_{\mathfrak B}+1)+\|\bfv_*\|_{W^{1,2}(W^{\frac32,2}(\partial\Omega))}^2+\|\bfg\|_{L^2(0,T)}\Big)\,.
\ea
\eeq{4.11}
In particular, the range of the map $M$ is contained in $\mathfrak B$. 
\EL{4.2}
{\em Proof.} By assumption, $(\bfphi,\bfpsi,\bfchi_\sub{\bfpsi})\in \mathfrak B$, and, by  \lemmref{4.1}, $\bfV\in \calw_\sharp^{2,\frac32}$. Thus, using  \eqref{WH2}, it follows at once $\bfF(\bfchi_\sub{\bfpsi},\bfxi_\sub{\bfphi} )\in L^2_\sharp(0,T)$ and that it satisfies \eqref{4.11}$_3$. The proof of \eqref{4.11}$_{1,2}$ becomes also quite straightforward, if we make use of the following continuous embedding properties
\be 
{\sf D}^{2,\frac32}\subset L^r(\Omega)\,,\ \ \mathcal W_\sharp^{2,\frac32}\subset L^q(L^\infty)\,,\  \ \mbox{for all $r\in [3,\infty]$ and $q\in [1,4)$\,.}
\eeq{4.11_1}
The second embedding is a particular case of \cite[Theorem 2.1]{Mallo}. To show the first one, it is enough to prove it for $r=\infty$. However, the latter follows from  the well-known inequality
\be
\|\bfphi\|_\infty\le c\,\big(\|\bfphi\|_3+|\bfphi|_{2,2}\big)\,;
\eeq{4.11_2}
see \cite[Theorem II.6.1 and Theorem II.9.1]{Gab}. Once the inequalities in \eqref{4.11} have been established, from \lemmref{Fm} and \lemmref{1.6_0} applied to \eqref{4.3} and \eqref{4.4}, respectively, we deduce that $(\bfu,\bfw,\bfchi_\sub{\bfu}):=M(\bfU)$ is in $\mathfrak B$, which completes the proof of the lemma. \par\hfill$\square$\par
\begin{thebibliography}{99}
\bibitem{Bab}Babenko, K.I., On the spectrum of a linearized problem on the flow of a viscous
incompressible fluid around a body (Russian). {\em Dokl. Akad. Nauk} SSSR, {\bf 262}, 64--68
(1982)
\bibitem{Bab1}Babenko, K.I., Periodic solutions of a problem of the flow of a viscous fluid around a
body. {\em Soviet Math. Dokl}. {\bf 25}, 211--216 (1982)
\bibitem{FaNe1}
Farwig, R., Neustupa, J.  On the spectrum of a Stokes-type
operator arising from flow around a rotating body. {\em Manuscripta
Math.} \textbf{122} (2007) 419--437
\bibitem{Gah} Galdi, G.P., On the motion of a rigid body in a viscous liquid: a mathematical analysis with applications. {\em Handbook of mathematical fluid dynamics}, Vol. I, 653--791, North-Holland, Amsterdam, 2002
\bibitem{Gafur}Galdi, G.P., Further properties of steady--state solutions to the Navier--Stokes problem past a three--dimensional obstacle, {\em J. Math. Phys.} {\bf 48}, 065207 (2007).
\bibitem{Gab}Galdi, G.P., {\em An introduction to the mathematical theory of the Navier--Stokes equations.
Steady-state problems}, Second edition. Springer Monographs in Mathematics,
Springer, New York (2011)
\bibitem{Gace}Galdi, G.P., Steady--state Navier--Stokes problem past a rotating body: geometric--functional properties and related questions. {\em Topics in mathematical fluid mechanics}, 109--197, Lecture Notes in Math., 2073, Fond. CIME/CIME Found. Subser., Springer, Heidelberg, 2013.
\bibitem{GaBif}Galdi, G.P., A time-periodic bifurcation theorem and its applications to Navier--Stokes flow past an obstacle, in {\em Mathematical Analysis of Viscous Incompressible Flow},
edited by T.~Hishida, R. I. M. S. Kokyuroku (Kyoto University, Japan, 2015), pp. 1--27
\bibitem{GaArma}Galdi, G.P., On bifurcating time-periodic flow of a Navier--Stokes liquid past a cylinder. {\em Arch. Ration. Mech. Anal.} {\bf 222} (2016), 285--315
\bibitem{GaJMP}Galdi, G.P., On the problem of steady bifurcation of a falling sphere in a Navier--Stokes liquid. {\em J. Math. Phys.} {\bf 61} (2020), no. 8, 083101, 13 pp
\bibitem{GaSP}Galdi G.P., On the self-propulsion of a rigid body in a viscous liquid by time-periodic boundary data. {\em J. Math. Fluid Mech.} {\bf 22} (2020), no. 4, Paper No. 61, 34 pp
\bibitem{GaMaH}
Galdi, G.P.,  Kyed, M., 
Time-periodic solutions to the Navier--Stokes equations.  {\em Handbook of mathematical analysis in mechanics of viscous fluids}, 509--578, Springer, Cham, 2018
\bibitem{GG}Gohberg, I., Goldberg, S. and Kaashoek, M.A., Classes of linear operators:
I. Operator Theory, Advances and Applications, Vol.49, Birkh\"auser
Verlag, Basel (1990)
\bibitem{HaIo}Haragus, M. and Iooss, G. {\em
Local bifurcations, center manifolds, and normal forms in infinite--dimensional dynamical systems.} 
Universitext. Springer--Verlag London, Ltd., London; EDP Sciences, Les Ulis, 2011
\bibitem{Hey} Heywood, J.G.,  
The Navier-Stokes equations: on the existence, regularity and decay of solutions, {\it Indiana Univ. Math. J.}, {\bf 29}, 639--681 (1980)
\bibitem{Jenny1}Jenny, M., Du\v{s}ek, J., Bouchet, G., Instabilities and transition of a sphere falling
or ascending freely in a Newtonian fluid. {\em J. Fluid Mech.} {\bf 508} (2004) 201--239
\bibitem{Kara}Karamanev, D., Chavarie, C., and  Mayer, R.,  Dynamics of the free rise of a light solid sphere
in liquid. {\em AIChE J.}, {\bf 42}, 1789--1792 (1996)
\bibitem{Nak}Nakamura, I., Steady wake behind a sphere, {\em Physics of Fluids} {\bf 19}  (1976) 5--8
\bibitem{Rao}Rao, M.A., {\em Rheology of Fluid and
Semisolid Foods
Principles and Applications}, Springer-Verlag, Second Edition (2014)
\bibitem{Satt}
Sattinger, D.H., {\em Topics in stability and bifurcation theory}. 
Lecture Notes in Mathematics, Vol. 309. Springer--Verlag, Berlin--New York, 1973
\bibitem{ST}Schweizer, T., A quick guide to better viscosity measurements of highly viscous fluids,
{\em Applied Rheology}, {\bf 14} (2004) 197--201
\bibitem{Scog}Scoggins, J.R., Aerodynamics of spherical balloon wind sensors, {\em J.
Geophys. Res.} {\bf 69}, 591-- 598 (1964)
\bibitem{ALS}Silvestre, A.L., On the self-propelled motion of a rigid body in a viscous liquid and on the attainability
of steady symmetric self-propelled motions, {\em J. Math. Fluid Mech}. {\bf 4} (2002) 285--326
\bibitem{Solo}Solonnikov, V.A., Estimates of the solutions of the nonstationary Navier--Stokes system.
{\em Zap. Naucn. Sem. Leningrad. Otdel. Mat. Inst. Steklov}. (LOMI) {\bf 38}, 153--201 (1973)
\bibitem{Tachi}Tachibana, M., The transient motion of a falling sphere in a
viscous fluid and the effects of side walls,{\em Mem.  Fac. Eng.
Fukui University}, {\bf 24}  (1976) 157--169 
\bibitem{Tan}Taneda, S., Kitasho, K., Experimental investigation of the wake behind a sphere at low Reynolds numbers. {\em J. Phys. Soc. Japan} {\bf 11}(1956) 1104--1108
\bibitem{Z1}Zeidler, E., {\em Nonlinear Functional Analysis and Applications}, Vol.4, Application to Mathematical Physics, Springer-Verlag, New York  (1988) 

\end{thebibliography}
